\documentclass[11pt,reqno]{amsart}
\calclayout

\usepackage{lmodern}

\usepackage{amssymb,mathtools}

\usepackage[hidelinks]{hyperref}

\newtheorem{Thm}{Theorem}[section]
\newtheorem{Lem}[Thm]{Lemma}

\newtheorem{Cor}[Thm]{Corollary}
\theoremstyle{remark}
\newtheorem{Rem}[Thm]{Remark}

\newcommand{\CC}{\mathbb{C}}
\newcommand{\ZZ}{\mathbb{Z}}

\newcommand{\Om}{\Omega}

\begin{document}

\title[Elliptic analogues of partial sums of very-well-poised ${}_8\psi_8$ series]{Contiguous relations for elliptic analogues of partial sums of very-well-poised $\boldsymbol{{}_8\psi_8}$ series}
\author{Takumi Maesaka}
\address{Faculty of Mathematics, Kyushu University, 744 Motooka, Nishi-ku, Fukuoka 819-0395, Japan}
\email{nozaki.takumi.912@s.kyushu-u.ac.jp}
\subjclass{33E05, 33D15, 33C20}
\keywords{Elliptic hypergeometric series, partial sums, contiguous relations, transformation formulas, Abel's lemma}

\begin{abstract}
We introduce finite elliptic analogues of partial sums of a balanced
very-well-poised ${}_8\psi_8$ series.  By combining Abel's lemma on summation by parts with the
Riemann theta relation, we establish two contiguous relations. Iterating the second relation gives two finite transformation formulas.
\end{abstract}

\maketitle

\section{Introduction}

Abel's lemma on summation by parts is a useful tool for deriving contiguous relations and transformations of hypergeometric series.
Its applications to ordinary and basic hypergeometric series were developed by Chu \cite{chu-classical,chu}, while partial $q$-series
transformations were studied systematically by Chu and Wang \cite{chu-wang}.
Closely related applications to nonterminating $q$-series and partial sums of Bailey's very-well-poised ${}_6\psi_6$ series were given by Chu and Zhang \cite{chu-zhang-nonterminating,chu-zhang-bailey}; see also the more recent general $q$-series transformations based on Abel's lemma due to Xu and Ma \cite{xu-ma}. For another application of Abel's lemma to contiguous relations for finite partial sums and their iteration, see Wang and Xu \cite{wang-xu}. Chu and Jia extended Abel's method to elliptic hypergeometric series \cite{chu-jia,chu-jia-elliptic}. Transformations of very-well-poised ${}_8\psi_8$ series go back to M. Jackson and Slater; see Gasper and Rahman \cite{gr} and Schlosser \cite{schlosser}.
In this paper, we use a similar approach to derive contiguous relations for elliptic analogues of partial sums of balanced very-well-poised ${}_8\psi_8$ series.
Specifically, we consider sums of the form
\[
\sum_{n=m}^{M-1}\frac{\theta(aq^{2n};p)(b,c,d,e,f,g;q,p)_n}{(aq/b,aq/c,aq/d,aq/e,aq/f,aq/g;q,p)_n}q^n, \qquad a^3q=bcdefg.
\]
Our first result is a contiguous relation that replaces $a$ by $aq$ and multiplies three of the six numerator parameters by $q$.
Its proof is a direct application of the finite form of Abel's lemma.
The second contiguous relation replaces $(f,g)$ by $(fq,g/q)$,
leaving $a,b,c,d,e$ unchanged.
By iterating this second relation, we obtain two finite transformation formulas for the elliptic partial sums.

The remainder of the paper is organized as follows.
Section \ref{sec:preliminaries} introduces the notation and states the finite form of Abel's lemma.
Section \ref{sec:contiguous} establishes two contiguous relations for the finite elliptic partial sums.
Section \ref{sec:iteration} derives two transformation formulas by iteration and obtains a finite partial-sum identity by comparing them.
Section \ref{sec:degeneration} presents the $p=0$ specializations and examines their nonterminating limits.
Section \ref{sec:classical-limit} gives the corresponding classical identities and analyzes the finite part in the boundary case.
The appendix presents well-poised analogues of Evans--Stanton-type formulas and their classical limits.

\section{Preliminaries}\label{sec:preliminaries}

Throughout the paper, let $0<|q|<1$ and $|p|<1$.  We write
\[
 (x;p)_\infty\coloneqq \prod_{j=0}^{\infty}(1-xp^j),
 \qquad
 (x_1,\ldots,x_r;p)_\infty\coloneqq \prod_{i=1}^r(x_i;p)_\infty
\]
and use the modified theta function
\[
 \theta(x;p)=(x,p/x;p)_\infty,
 \qquad
 \theta(x_1,\ldots,x_r;p)=\prod_{j=1}^r\theta(x_j;p).
\]
It satisfies
\begin{equation}\label{eq:theta-inversion}
 \theta(px;p)=\theta(x^{-1};p)=-x^{-1}\theta(x;p).
\end{equation}
We shall also use the Riemann theta relation \cite[Eq.~(3.1)]{warnaar}
\begin{align}
&\theta(xy,x/y,uv,u/v;p)-\theta(xv,x/v,uy,u/y;p)=\frac{u}{y}\theta(yv,y/v,xu,x/u;p).
\label{eq:riemann}
\end{align}
For standard background on theta functions and elliptic hypergeometric series, see Warnaar \cite{warnaar} and Gasper and Rahman \cite{gr}.

For $x\in\CC^\times$ and $n\in\ZZ$, define
\[
 (x;q,p)_n=
 \begin{cases}
 \displaystyle\prod_{j=0}^{n-1}\theta(xq^j;p),&n>0,\\
 1,&n=0,\\
 \displaystyle\prod_{j=1}^{-n}\theta(xq^{-j};p)^{-1},&n<0,
 \end{cases}
\]
and abbreviate
\[
 (x_1,\ldots,x_r;q,p)_n=\prod_{j=1}^r(x_j;q,p)_n.
\]
The shifted factorials satisfy the elementary product relation
\begin{equation}\label{eq:shifted-factorial-product}
 (x;q,p)_{r+s}=(x;q,p)_r(xq^r;q,p)_s
 \qquad(r,s\in\ZZ)
\end{equation}
We shall use \eqref{eq:shifted-factorial-product} repeatedly.

Let $m,M\in\ZZ$ with $m\leq M$.  Under the balancing condition
\begin{equation}\label{eq:balance}
 a^3q=bcdefg,
\end{equation}
we define the finite elliptic partial sum
\begin{align}
&\Om_m^M(a\mid b,c,d,e,f,g;q,p)\coloneqq\sum_{n=m}^{M-1}\frac{\theta(aq^{2n};p)(b,c,d,e,f,g;q,p)_n}{(aq/b,aq/c,aq/d,aq/e,aq/f,aq/g;q,p)_n}q^n.\label{eq:omega-def}
\end{align}
Thus $\Om_M^M(a\mid b,c,d,e,f,g;q,p)=0$ by the empty-sum convention. All parameters are assumed to be generic, so that the denominators occurring in the formulas do not vanish. Since the sums in Sections \ref{sec:preliminaries}--\ref{sec:iteration} are finite, the identities in those sections extend meromorphically in the parameters. The identities involving infinite series considered later are subject to the stated convergence conditions.

We record the finite form of Abel's lemma used below; see also \cite{chu-classical, chu}.

\begin{Lem}\label{lem:abel}
Let $m\leq M$ be integers.
For any sequences $(A_n)$ and $(B_n)$ for which
the following expressions are defined, we have
\begin{align}
\sum_{n=m}^{M-1}B_n(A_n-A_{n-1})&=A_{M-1}B_M-A_{m-1}B_m+\sum_{n=m}^{M-1}A_n(B_n-B_{n+1}).
\label{eq:abel}
\end{align}
\end{Lem}

\section{Two finite elliptic contiguous relations}\label{sec:contiguous}

Chu and Jia applied a modified form of Abel's lemma to terminating theta hypergeometric series \cite{chu-jia}. We use the same general approach, combining Abel's lemma with the Riemann theta relation, to obtain a contiguous relation for sums over an arbitrary finite interval $m\leq n<M$.

We replace $(e,f,g)$ by $(eq,fq,gq)$.  For $r\in\ZZ$, put
\begin{align*}
E_r^{(q,p)}(a\mid b,c,d,e,f,g)
&\coloneqq
\frac{efg\,\theta(efgq^r/a;p)}
{a^2\theta(ef/a,eg/a,fg/a;p)}\notag\\
&\quad\cdot
\frac{(b,c,d,e,f,g;q,p)_r}
{(aq/e,aq/f,aq/g;q,p)_{r-1}
 (aq/b,aq/c,aq/d;q,p)_r},
\end{align*}
and
\begin{equation}\label{eq:R-def}
R^{(q,p)}\coloneqq
\frac{\theta(e,f,g,aq/bc,aq/bd,aq/cd;p)}
{\theta(ef/a,eg/a,fg/a,aq/b,aq/c,aq/d;p)}.
\end{equation}

\begin{Thm}\label{thm:first-contiguous}
Let $m\leq M$ be integers.  Under \eqref{eq:balance},
\begin{align}
\Om_m^M(a\mid b,c,d,e,f,g;q,p)&=E_m^{(q,p)}(a\mid b,c,d,e,f,g)-E_M^{(q,p)}(a\mid b,c,d,e,f,g)\notag\\
&\quad+R^{(q,p)}\Om_m^M(aq\mid b,c,d,eq,fq,gq;q,p).
\label{eq:first-contiguous}
\end{align}
\end{Thm}

\begin{proof}
Define
\begin{align*}
A_n&=\frac{(eq,fq,gq,a^2q/efg;q,p)_n}{(aq/e,aq/f,aq/g,efgq/a;q,p)_n},\\
B_n&=\frac{(b,c,d,a^2q^2/bcd;q,p)_n}{(aq/b,aq/c,aq/d,bcd/aq;q,p)_n}.
\end{align*}
The two quotients to which we apply \eqref{eq:riemann} are
\begin{align*}
 \frac{A_n}{A_{n-1}}&=\frac{\theta(eq^n,fq^n,gq^n,a^2q^n/efg;p)}{\theta(aq^n/e,aq^n/f,aq^n/g,efgq^n/a;p)},\\
 \frac{B_{n+1}}{B_n}&=\frac{\theta(bq^n,cq^n,dq^n,a^2q^{n+2}/bcd;p)}{\theta(aq^{n+1}/b,aq^{n+1}/c,aq^{n+1}/d,bcdq^{n-1}/a;p)}.
\end{align*}
Expressing each difference over a common denominator and applying \eqref{eq:riemann} gives
\begin{align}
A_n-A_{n-1}&=C_A\frac{\theta(aq^{2n};p)(e,f,g,a^2/efg;q,p)_n}{(aq/e,aq/f,aq/g,efgq/a;q,p)_n}q^n,
\notag\\
B_n-B_{n+1}&=C_B\frac{\theta(aq^{2n+1};p)(b,c,d,a^2q^2/bcd;q,p)_n}{(aq^2/b,aq^2/c,aq^2/d,bcd/a;q,p)_n}q^n,
\notag
\end{align}
where
\begin{align}
C_A&=-\frac{a^2\theta(ef/a,eg/a,fg/a;p)}{efg\,\theta(e,f,g,a^2/efg;p)},\notag\\
C_B&=\frac{\theta(aq/bc,aq/bd,aq/cd;p)}{\theta(aq/b,aq/c,aq/d,aq/bcd;p)}.\notag
\end{align}
Using the balancing condition together with
\eqref{eq:shifted-factorial-product}, we obtain
\begin{align*}
B_n(A_n-A_{n-1})&=C_A\frac{\theta(aq^{2n};p)(b,c,d,e,f,g;q,p)_n}{(aq/b,aq/c,aq/d,aq/e,aq/f,aq/g;q,p)_n}q^n,\\
A_n(B_n-B_{n+1})&=C_B\frac{\theta(aq^{2n+1};p)(b,c,d,eq,fq,gq;q,p)_n}{(aq^2/b,aq^2/c,aq^2/d,aq/e,aq/f,aq/g;q,p)_n}q^n.
\end{align*}
Moreover, by \eqref{eq:theta-inversion},
\begin{align*}
 \frac{C_B}{C_A}=R^{(q,p)},\qquad-\frac{A_{r-1}B_r}{C_A}=E_r^{(q,p)}(a\mid b,c,d,e,f,g).
\end{align*}
Substitution of these identities into \eqref{eq:abel}, followed by division by $C_A$, gives \eqref{eq:first-contiguous}.
\end{proof}

The second relation preserves $a,b,c,d,e$ and moves the pair $(f,g)$ in opposite directions. Define
\begin{align}
F_r^{(q,p)}(a\mid b,c,d,e,f,g)&\coloneqq\frac{\theta(g/fq;p)}{f\theta(a/bf,a/cf,a/df,a/ef;p)}\notag\\
&\quad\cdot\frac{(b,c,d,e,f;q,p)_r(g;q,p)_{r-1}}{(aq/b,aq/c,aq/d,aq/e,aq/f;q,p)_{r-1}(aq/g;q,p)_r},\notag
\end{align}
and
\begin{equation}\label{eq:K-def}
K^{(q,p)}\coloneqq\frac{g^2}{f^2q^2}\frac{\theta(f,a/f,aq/bg,aq/cg,aq/dg,aq/eg;p)}{\theta(g/q,aq/g,a/bf,a/cf,a/df,a/ef;p)}.
\end{equation}

\begin{Thm}\label{thm:second-contiguous}
Let $m\leq M$ be integers.  Under \eqref{eq:balance},
\begin{align}
\Om_m^M(a\mid b,c,d,e,f,g;q,p)&=F_m^{(q,p)}(a\mid b,c,d,e,f,g)-F_M^{(q,p)}(a\mid b,c,d,e,f,g)\notag\\
&\quad+K^{(q,p)}\Om_m^M(a\mid b,c,d,e,fq,g/q;q,p).\label{eq:second-contiguous}
\end{align}
\end{Thm}

\begin{proof}
We compare two permuted instances of Theorem~\ref{thm:first-contiguous}.  To make the comparison transparent, choose a
shifted triple $(x,y,z)$ from the six parameters and let $(u,v,w)$ be the complementary triple.  We write
\begin{align*}
 E_r^{(q,p);x,y,z}&\coloneqq E_r^{(q,p)}(a\mid u,v,w,x,y,z),
\end{align*}
and let $R^{(q,p);x,y,z}$ denote the expression in \eqref{eq:R-def} with $(b,c,d,e,f,g)=(u,v,w,x,y,z)$. The order within either triple is immaterial.
In particular, the complementary triples associated with $(d,e,g)$ and $(d,e,f)$ are $(b,c,f)$ and $(b,c,g)$, respectively.  Symmetry of
\eqref{eq:omega-def} gives
\begin{align*}
\Om_m^M(a\mid b,c,d,e,f,g;q,p)&=E_m^{(q,p);d,e,g}-E_M^{(q,p);d,e,g}+R^{(q,p);d,e,g}\Om_m^M(aq\mid b,c,dq,eq,f,gq;q,p),\\
\Om_m^M(a\mid b,c,d,e,f,g;q,p)&=E_m^{(q,p);d,e,f}-E_M^{(q,p);d,e,f}+R^{(q,p);d,e,f}\Om_m^M(aq\mid b,c,dq,eq,fq,g;q,p).
\end{align*}
Subtracting these identities yields
\begin{align}
&R^{(q,p);d,e,g}\Om_m^M(aq\mid b,c,dq,eq,f,gq;q,p)-R^{(q,p);d,e,f}\Om_m^M(aq\mid b,c,dq,eq,fq,g;q,p)\notag\\
&\quad=(E_m^{(q,p);d,e,f}-E_m^{(q,p);d,e,g})-(E_M^{(q,p);d,e,f}-E_M^{(q,p);d,e,g}).\label{eq:comparison}
\end{align}
In \eqref{eq:comparison}, first replace $(a,d,e)$ by $(a/q,d/q,e/q)$ and then replace $g$ by $g/q$.
The resulting balancing condition is \eqref{eq:balance}. After these substitutions, we have
\[
 \frac{R^{(q,p);d,e,f}}{R^{(q,p);d,e,g}}=K^{(q,p)},\qquad\frac{E_r^{(q,p);d,e,f}-E_r^{(q,p);d,e,g}}{R^{(q,p);d,e,g}}=F_r^{(q,p)}(a\mid b,c,d,e,f,g).
\]
Dividing the resulting identity by $R^{(q,p);d,e,g}$ therefore gives \eqref{eq:second-contiguous}.
\end{proof}

\section{Iterated transformations}\label{sec:iteration}

The method of deriving transformation formulas by iterating contiguous relations has been used in \cite{chu-jia,chu-jia-elliptic,chu-wang, chu-zhang-nonterminating,chu-zhang-bailey}. We now apply this method to Theorem~\ref{thm:second-contiguous}. For a nonnegative integer $N$, define
\begin{align*}
W_N^{(q,p)}\coloneqq\frac{(f,f/a,aq/bg,aq/cg,aq/dg,aq/eg;q,p)_N}{(q/g,aq/g,bf/a,cf/a,df/a,ef/a;q,p)_N}.
\end{align*}

\begin{Thm}\label{thm:fixed-lower}
Let $m\leq M$ be integers.  Under \eqref{eq:balance},
\begin{align}
&W_N^{(q,p)}\Om_m^M(a\mid b,c,d,e,fq^N,gq^{-N};q,p)\notag\\
&=\Om_m^M(a\mid b,c,d,e,f,g;q,p)\notag\\
&\quad+\frac{f}{a}\frac{(b,c,d,e,f;q,p)_m(g;q,p)_{m-1}}{(aq/b,aq/c,aq/d,aq/e,aq/f;q,p)_{m-1}(aq/g;q,p)_m}\notag\\
&\qquad\cdot\sum_{n=0}^{N-1}\frac{\theta(fq^{2n+1}/g;p)(aq/bg,aq/cg,aq/dg,aq/eg,fq^{1-m}/a,fq^m;q,p)_n}{(bf/a,cf/a,df/a,ef/a;q,p)_{n+1}(q^{2-m}/g,aq^{m+1}/g;q,p)_n}q^n\notag\\
&\quad-\frac{f}{a}\frac{(b,c,d,e,f;q,p)_M(g;q,p)_{M-1}}{(aq/b,aq/c,aq/d,aq/e,aq/f;q,p)_{M-1}(aq/g;q,p)_M}\notag\\
&\qquad\cdot\sum_{n=0}^{N-1}\frac{\theta(fq^{2n+1}/g;p)(aq/bg,aq/cg,aq/dg,aq/eg,fq^{1-M}/a,fq^M;q,p)_n}{(bf/a,cf/a,df/a,ef/a;q,p)_{n+1}(q^{2-M}/g,aq^{M+1}/g;q,p)_n}q^n.\label{eq:fixed-lower}
\end{align}
\end{Thm}

\begin{proof}
Let $K_n^{(q,p)}$ denote the expression in \eqref{eq:K-def} with $(f,g)$ replaced by $(fq^n,gq^{-n})$.  Direct use of \eqref{eq:shifted-factorial-product} gives $W_{n+1}^{(q,p)}=W_n^{(q,p)}K_n^{(q,p)}$.  Moreover, for every $r\in\ZZ$,
\begin{align}
&-W_n^{(q,p)}F_r^{(q,p)}(a\mid b,c,d,e,fq^n,gq^{-n})\notag\\
&\quad=\frac{f}{a}\frac{(b,c,d,e,f;q,p)_r(g;q,p)_{r-1}}{(aq/b,aq/c,aq/d,aq/e,aq/f;q,p)_{r-1}(aq/g;q,p)_r}\notag\\
&\qquad\cdot\frac{\theta(fq^{2n+1}/g;p)(aq/bg,aq/cg,aq/dg,aq/eg,fq^{1-r}/a,fq^r;q,p)_n}{(bf/a,cf/a,df/a,ef/a;q,p)_{n+1}(q^{2-r}/g,aq^{r+1}/g;q,p)_n}q^n.\label{eq:weighted-boundary}
\end{align}
Apply \eqref{eq:second-contiguous} with $(f,g)=(fq^n,gq^{-n})$ and multiply by
$W_n^{(q,p)}$.  If
\[
 S_n=W_n^{(q,p)}\Om_m^M(a\mid b,c,d,e,fq^n,gq^{-n};q,p),
\]
then $W_{n+1}^{(q,p)}=W_n^{(q,p)}K_n^{(q,p)}$ gives
\[
 S_n-S_{n+1}=W_n^{(q,p)}F_m^{(q,p)}(a\mid b,c,d,e,fq^n,gq^{-n})-W_n^{(q,p)}F_M^{(q,p)}(a\mid b,c,d,e,fq^n,gq^{-n}).
\]
Summing this equality for $0\leq n<N$ and using
\eqref{eq:weighted-boundary} with $r=m,M$ proves
\eqref{eq:fixed-lower}.
\end{proof}

\begin{Thm}\label{thm:moving-lower}
Let $N$ and $M$ be integers with $0\leq N\leq M$.  Under
\eqref{eq:balance},
\begin{align}
&W_N^{(q,p)}\Om_N^M(a\mid b,c,d,e,fq^N,gq^{-N};q,p)\notag\\
&=\Om_0^M(a\mid b,c,d,e,f,g;q,p)\notag\\
&\quad-\frac{f}{a}\frac{(b,c,d,e,f;q,p)_M(g;q,p)_{M-1}}{(aq/b,aq/c,aq/d,aq/e,aq/f;q,p)_{M-1}(aq/g;q,p)_M}\notag\\
&\qquad\cdot\sum_{n=0}^{N-1}\frac{\theta(fq^{2n+1}/g;p)(aq/bg,aq/cg,aq/dg,aq/eg,fq^{1-M}/a,fq^M;q,p)_n}{(bf/a,cf/a,df/a,ef/a;q,p)_{n+1}(q^{2-M}/g,aq^{M+1}/g;q,p)_n}q^n\notag\\
&\quad+\frac{f}{a}\sum_{n=0}^{N-1}\frac{\theta(fq^{2n+1}/g;p)(b,c,d,e;q,p)_{n+1}}{(aq/b,aq/c,aq/d,aq/e;q,p)_n}\notag\\
&\qquad\cdot\frac{(aq/bg,aq/cg,aq/dg,aq/eg;q,p)_n}{(bf/a,cf/a,df/a,ef/a;q,p)_{n+1}}\frac{(f;q,p)_{2n+1}}{(aq/g;q,p)_{2n+1}}\left(\frac{fg}{a}\right)^n\notag\\
&\quad-\sum_{n=0}^{N-1}\frac{\theta(aq^{2n};p)(b,c,d,e,aq/bg,aq/cg,aq/dg,aq/eg;q,p)_n}{(aq/b,aq/c,aq/d,aq/e,bf/a,cf/a,df/a,ef/a;q,p)_n}\frac{(f;q,p)_{2n}}{(aq/g;q,p)_{2n}}\left(\frac{fg}{a}\right)^n.\label{eq:moving-lower}
\end{align}
\end{Thm}

\begin{proof}
Put
\[
 T_n=W_n^{(q,p)}\Om_n^M(a\mid b,c,d,e,fq^n,gq^{-n};q,p).
\]
Separate the $n$th summand from $\Om_n^M(a\mid b,c,d,e,fq^n,gq^{-n};q,p)$ and apply Theorem~\ref{thm:second-contiguous} to the remaining sum $\Om_{n+1}^M(a\mid b,c,d,e,fq^n,gq^{-n};q,p)$. Putting $r=n+1$ in \eqref{eq:weighted-boundary} and simplifying by \eqref{eq:shifted-factorial-product}, we obtain
\begin{align}
&-W_n^{(q,p)}F_{n+1}^{(q,p)}(a\mid b,c,d,e,fq^n,gq^{-n})\notag\\
&\quad=\frac{f}{a}\frac{\theta(fq^{2n+1}/g;p)(b,c,d,e;q,p)_{n+1}}{(aq/b,aq/c,aq/d,aq/e;q,p)_n}\notag\\
&\qquad\cdot\frac{(aq/bg,aq/cg,aq/dg,aq/eg;q,p)_n}{(bf/a,cf/a,df/a,ef/a;q,p)_{n+1}}\frac{(f;q,p)_{2n+1}}{(aq/g;q,p)_{2n+1}}\left(\frac{fg}{a}\right)^n.\label{eq:moving-boundary}
\end{align}
The weighted $n$th summand is
\begin{align}
&W_n^{(q,p)}\bigl[\text{$n$th summand of }\Om_n^M(a\mid b,c,d,e,fq^n,gq^{-n};q,p)\bigr]\notag\\
&\quad=\frac{\theta(aq^{2n};p)(b,c,d,e,aq/bg,aq/cg,aq/dg,aq/eg;q,p)_n}{(aq/b,aq/c,aq/d,aq/e,bf/a,cf/a,df/a,ef/a;q,p)_n}\frac{(f;q,p)_{2n}}{(aq/g;q,p)_{2n}}\left(\frac{fg}{a}\right)^n.
\label{eq:weighted-summand}
\end{align}
Since $W_{n+1}^{(q,p)}=W_n^{(q,p)}K_n^{(q,p)}$, Theorem~\ref{thm:second-contiguous} now gives
\begin{align*}
T_n-T_{n+1}&=-W_n^{(q,p)}F_M^{(q,p)}(a\mid b,c,d,e,fq^n,gq^{-n})\\
&\quad+W_n^{(q,p)}F_{n+1}^{(q,p)}(a\mid b,c,d,e,fq^n,gq^{-n})\\
&\quad+W_n^{(q,p)}\bigl[\text{$n$th summand of }\Om_n^M(a\mid b,c,d,e,fq^n,gq^{-n};q,p)\bigr].
\end{align*}
Summing over $0\leq n<N$ and substituting
\eqref{eq:weighted-boundary} with $r=M$, together with
\eqref{eq:moving-boundary} and \eqref{eq:weighted-summand}, proves
\eqref{eq:moving-lower}.
\end{proof}

The two iterated formulas immediately imply a closed relation among finite
sums. 

\begin{Cor}\label{cor:finite-sum}
For every nonnegative integer $N$, under \eqref{eq:balance},
\begin{align}
&W_N^{(q,p)}\Om_0^N(a\mid b,c,d,e,fq^N,gq^{-N};q,p)\notag\\
&\quad=\frac{f}{a}\frac{\theta(a/b,a/c,a/d,a/e,a/f;p)}{\theta(g/q;p)}\notag\\
&\qquad\cdot\sum_{n=0}^{N-1}\frac{\theta(fq^{2n+1}/g;p)(aq/bg,aq/cg,aq/dg,aq/eg,fq/a,f;q,p)_n}{(bf/a,cf/a,df/a,ef/a;q,p)_{n+1}(q^2/g,aq/g;q,p)_n}q^n\notag\\
&\quad-\frac{f}{a}\sum_{n=0}^{N-1}\frac{\theta(fq^{2n+1}/g;p)(b,c,d,e;q,p)_{n+1}}{(aq/b,aq/c,aq/d,aq/e;q,p)_n}\notag\\
&\qquad\cdot\frac{(aq/bg,aq/cg,aq/dg,aq/eg;q,p)_n}{(bf/a,cf/a,df/a,ef/a;q,p)_{n+1}}\frac{(f;q,p)_{2n+1}}{(aq/g;q,p)_{2n+1}}\left(\frac{fg}{a}\right)^n\notag\\
&\quad+\sum_{n=0}^{N-1}\frac{\theta(aq^{2n};p)(b,c,d,e,aq/bg,aq/cg,aq/dg,aq/eg;q,p)_n}{(aq/b,aq/c,aq/d,aq/e,bf/a,cf/a,df/a,ef/a;q,p)_n}\frac{(f;q,p)_{2n}}{(aq/g;q,p)_{2n}}\left(\frac{fg}{a}\right)^n.\label{eq:finite-sum}
\end{align}
\end{Cor}

\begin{proof}
Apply Theorem~\ref{thm:fixed-lower} with $(m,M)=(0,N)$.  Next apply Theorem~\ref{thm:moving-lower} with $M=N$; its left-hand side then vanishes by the empty-sum convention.  Eliminate the common untransformed partial sum from the two identities. Simplifying the term corresponding to $m=0$ using the negative-index convention for elliptic shifted factorials then gives \eqref{eq:finite-sum}.
\end{proof}

\section{The \texorpdfstring{\(p=0\)}{p=0} specialization}
\label{sec:degeneration}

At $p=0$ one has
\[
 \theta(x;0)=1-x,\qquad (x;q,0)_n=(x;q)_n.
\]
For integers $m\leq M$, we therefore use the notation
\begin{align}
\Om_m^M(a\mid b,c,d,e,f,g;q)&\coloneqq \Om_m^M(a\mid b,c,d,e,f,g;q,0)\notag\\
&=\sum_{n=m}^{M-1}\frac{(1-aq^{2n})(b,c,d,e,f,g;q)_n}{(aq/b,aq/c,aq/d,aq/e,aq/f,aq/g;q)_n}q^n.\notag
\end{align}
To state the $p=0$ specialization of Theorem~\ref{thm:first-contiguous}, define
\begin{align}
E_r^{(q)}&\coloneqq\frac{efg(1-efgq^r/a)}{a^2(1-ef/a)(1-eg/a)(1-fg/a)}\notag\\
&\quad\cdot\frac{(b,c,d,e,f,g;q)_r}{(aq/e,aq/f,aq/g;q)_{r-1}(aq/b,aq/c,aq/d;q)_r}.\notag
\end{align}
\begin{Cor}[$p=0$ specialization of Theorem~\ref{thm:first-contiguous}]\label{cor:first-contiguous-q}
Let $m\leq M$ be integers.  Under \eqref{eq:balance},
\begin{align}
\Om_m^M(a\mid b,c,d,e,f,g;q)&=E_m^{(q)}-E_M^{(q)}\notag\\
&\quad+\frac{(1-e)(1-f)(1-g)(1-aq/bc)(1-aq/bd)(1-aq/cd)}{(1-ef/a)(1-eg/a)(1-fg/a)(1-aq/b)(1-aq/c)(1-aq/d)}\notag\\
&\qquad\cdot\Om_m^M(aq\mid b,c,d,eq,fq,gq;q).\notag
\end{align}
\end{Cor}

Likewise, define
\begin{align}
F_r^{(q)}&\coloneqq\frac{1-g/fq}{f(1-a/bf)(1-a/cf)(1-a/df)(1-a/ef)}\notag\\
&\quad\cdot\frac{(b,c,d,e,f;q)_r(g;q)_{r-1}}{(aq/b,aq/c,aq/d,aq/e,aq/f;q)_{r-1}(aq/g;q)_r}.
\notag
\end{align}
\begin{Cor}[$p=0$ specialization of Theorem
\ref{thm:second-contiguous}]\label{cor:second-contiguous-q}
Let $m\leq M$ be integers.  Under \eqref{eq:balance},
\begin{align*}
&\Om_m^M(a\mid b,c,d,e,f,g;q)\\
&=F_m^{(q)}-F_M^{(q)}\notag\\
&\quad+\frac{g^2}{f^2q^2}\frac{(1-f)(1-a/f)(1-aq/bg)(1-aq/cg)(1-aq/dg)(1-aq/eg)}{(1-g/q)(1-aq/g)(1-a/bf)(1-a/cf)(1-a/df)(1-a/ef)}\notag\\
&\qquad\cdot\Om_m^M(a\mid b,c,d,e,fq,g/q;q).
\end{align*}
\end{Cor}

For convenience, set
\begin{equation*}
W_N^{(q)}\coloneqq\frac{(f,f/a,aq/bg,aq/cg,aq/dg,aq/eg;q)_N}{(q/g,aq/g,bf/a,cf/a,df/a,ef/a;q)_N}.
\end{equation*}
We also write
\begin{equation*}
W_\infty^{(q)}\coloneqq\frac{(f,f/a,aq/bg,aq/cg,aq/dg,aq/eg;q)_\infty}{(q/g,aq/g,bf/a,cf/a,df/a,ef/a;q)_\infty}.
\end{equation*}
Here $(x;q)_\infty\coloneqq \prod_{j=0}^{\infty}(1-xq^j)$ and multiple arguments denote products, as usual. The $p=0$ specializations of the three results in Section \ref{sec:iteration} are given explicitly below.

These identities should be compared with the well-poised partial-$q$-series transformations of Chu and Wang \cite{chu-wang} and with the transformations of partial sums of Bailey's ${}_6\psi_6$ series obtained by Chu and Zhang \cite{chu-zhang-bailey}.  In contrast with the complete ${}_8\psi_8$ transformations in \cite[Chapter~5]{gr} and \cite{schlosser}, the finite formulas below retain the lower and upper boundary terms.

\begin{Cor}[$p=0$ specialization of Theorem~\ref{thm:fixed-lower}]
\label{cor:fixed-lower-q}
Let $N$ be a nonnegative integer and let $m\leq M$ be integers.  Under \eqref{eq:balance},
\begin{align}
&W_N^{(q)}\Om_m^M(a\mid b,c,d,e,fq^N,gq^{-N};q)\notag\\
&=\Om_m^M(a\mid b,c,d,e,f,g;q)\notag\\
&\quad+
\frac{f}{a}\frac{(b,c,d,e,f;q)_m(g;q)_{m-1}}{(aq/b,aq/c,aq/d,aq/e,aq/f;q)_{m-1}(aq/g;q)_m}\notag\\
&\qquad\cdot\sum_{n=0}^{N-1}\frac{(1-fq^{2n+1}/g)(aq/bg,aq/cg,aq/dg,aq/eg,fq^{1-m}/a,fq^m;q)_n}{(bf/a,cf/a,df/a,ef/a;q)_{n+1}(q^{2-m}/g,aq^{m+1}/g;q)_n}q^n\notag\\
&\quad-\frac{f}{a}\frac{(b,c,d,e,f;q)_M(g;q)_{M-1}}{(aq/b,aq/c,aq/d,aq/e,aq/f;q)_{M-1}(aq/g;q)_M}\notag\\
&\qquad\cdot\sum_{n=0}^{N-1}\frac{(1-fq^{2n+1}/g)(aq/bg,aq/cg,aq/dg,aq/eg,fq^{1-M}/a,fq^M;q)_n}{(bf/a,cf/a,df/a,ef/a;q)_{n+1}(q^{2-M}/g,aq^{M+1}/g;q)_n}q^n.\label{eq:fixed-lower-q}
\end{align}
\end{Cor}

For later use, define
\begin{align*}
\Om_m^\infty(a\mid b,c,d,e,f,g;q)&\coloneqq\sum_{k=m}^{\infty}\frac{(1-aq^{2k})(b,c,d,e,f,g;q)_k}{(aq/b,aq/c,aq/d,aq/e,aq/f,aq/g;q)_k}q^k.
\end{align*}
Setting $M=N$ in \eqref{eq:fixed-lower-q} and letting $N\to\infty$, we obtain the following corollary.
\begin{Cor}\label{cor:fixed-lower-q-infinity}
Let $m\in\ZZ$, and assume, in addition to \eqref{eq:balance}, that $|fg/a|<1$.  Then
\begin{align*}
&W_\infty^{(q)}\sum_{k=m}^{\infty}\frac{(1-aq^{2k})(b,c,d,e;q)_k}{(aq/b,aq/c,aq/d,aq/e;q)_k}\left(\frac{fg}{a}\right)^k\notag\\
&\quad=\Om_m^\infty(a\mid b,c,d,e,f,g;q)\notag\\
&\qquad+\frac{f}{a}\frac{(b,c,d,e,f;q)_m(g;q)_{m-1}}{(aq/b,aq/c,aq/d,aq/e,aq/f;q)_{m-1}(aq/g;q)_m}\notag\\
&\qquad\quad\cdot\sum_{n=0}^{\infty}\frac{(1-fq^{2n+1}/g)(aq/bg,aq/cg,aq/dg,aq/eg,fq^{1-m}/a,fq^m;q)_n}{(bf/a,cf/a,df/a,ef/a;q)_{n+1}(q^{2-m}/g,aq^{m+1}/g;q)_n}q^n\notag\\
&\qquad-\frac{f}{a}\frac{(b,c,d,e,f,g;q)_\infty}{(aq/b,aq/c,aq/d,aq/e,aq/f,aq/g;q)_\infty}\notag\\
&\qquad\quad\cdot\sum_{n=0}^{\infty}\frac{(1-fq^{2n+1}/g)(aq/bg,aq/cg,aq/dg,aq/eg;q)_n}{(bf/a,cf/a,df/a,ef/a;q)_{n+1}}\left(\frac{fg}{a}\right)^n.
\end{align*}
\end{Cor}

\begin{proof}
Set $M=N$ in \eqref{eq:fixed-lower-q}.  For every fixed $k\in\ZZ$,
\begin{equation*}
 \lim_{N\to\infty}q^k\frac{(fq^N,gq^{-N};q)_k}{(aq^{1-N}/f,aq^{N+1}/g;q)_k}=\left(\frac{fg}{a}\right)^k,
\end{equation*}
and $W_N^{(q)}\to W_\infty^{(q)}$. For each fixed $n\geq0$, we also have
\begin{equation*}
 \lim_{N\to\infty}q^n\frac{(fq^{1-N}/a,fq^N;q)_n}{(q^{2-N}/g,aq^{N+1}/g;q)_n}=\left(\frac{fg}{a}\right)^n.
\end{equation*}
The remaining prefactor tends to
\[
 \frac{f}{a}\frac{(b,c,d,e,f,g;q)_\infty}{(aq/b,aq/c,aq/d,aq/e,aq/f,aq/g;q)_\infty}.
\]
Substitution of these limits into \eqref{eq:fixed-lower-q} proves the claim.
\end{proof}

\begin{Cor}[$p=0$ specialization of Theorem~\ref{thm:moving-lower}]
\label{cor:moving-lower-q}
Let $N$ and $M$ be integers with $0\leq N\leq M$.  Under
\eqref{eq:balance},
\begin{align*}
&W_N^{(q)}\Om_N^M(a\mid b,c,d,e,fq^N,gq^{-N};q)\\
&=\Om_0^M(a\mid b,c,d,e,f,g;q)\\
&\quad-\frac{f}{a}\frac{(b,c,d,e,f;q)_M(g;q)_{M-1}}{(aq/b,aq/c,aq/d,aq/e,aq/f;q)_{M-1}(aq/g;q)_M}\\
&\qquad\cdot\sum_{n=0}^{N-1}\frac{(1-fq^{2n+1}/g)(aq/bg,aq/cg,aq/dg,aq/eg,fq^{1-M}/a,fq^M;q)_n}{(bf/a,cf/a,df/a,ef/a;q)_{n+1}(q^{2-M}/g,aq^{M+1}/g;q)_n}q^n\\
&\quad+\frac{f}{a}\sum_{n=0}^{N-1}\frac{(1-fq^{2n+1}/g)(b,c,d,e;q)_{n+1}}
{(aq/b,aq/c,aq/d,aq/e;q)_n}\\
&\qquad\cdot\frac{(aq/bg,aq/cg,aq/dg,aq/eg;q)_n}{(bf/a,cf/a,df/a,ef/a;q)_{n+1}}\frac{(f;q)_{2n+1}}{(aq/g;q)_{2n+1}}\left(\frac{fg}{a}\right)^n\\
&\quad-\sum_{n=0}^{N-1}\frac{(1-aq^{2n})(b,c,d,e,aq/bg,aq/cg,aq/dg,aq/eg;q)_n}{(aq/b,aq/c,aq/d,aq/e,bf/a,cf/a,df/a,ef/a;q)_n}\frac{(f;q)_{2n}}{(aq/g;q)_{2n}}\left(\frac{fg}{a}\right)^n.
\end{align*}
\end{Cor}

\begin{Cor}[$p=0$ specialization of Corollary \ref{cor:finite-sum}]
\label{cor:finite-sum-q}
For every nonnegative integer $N$, under \eqref{eq:balance},
\begin{align}
&W_N^{(q)}\Om_0^N(a\mid b,c,d,e,fq^N,gq^{-N};q)\notag\\
&\quad=\frac{f}{a}\frac{(1-a/b)(1-a/c)(1-a/d)(1-a/e)(1-a/f)}{1-g/q}\notag\\
&\qquad\cdot\sum_{n=0}^{N-1}\frac{(1-fq^{2n+1}/g)(aq/bg,aq/cg,aq/dg,aq/eg,fq/a,f;q)_n}{(bf/a,cf/a,df/a,ef/a;q)_{n+1}(q^2/g,aq/g;q)_n}q^n\notag\\
&\quad-\frac{f}{a}\sum_{n=0}^{N-1}\frac{(1-fq^{2n+1}/g)(b,c,d,e;q)_{n+1}}{(aq/b,aq/c,aq/d,aq/e;q)_n}\notag\\
&\qquad\cdot\frac{(aq/bg,aq/cg,aq/dg,aq/eg;q)_n}{(bf/a,cf/a,df/a,ef/a;q)_{n+1}}\frac{(f;q)_{2n+1}}{(aq/g;q)_{2n+1}}\left(\frac{fg}{a}\right)^n\notag\\
&\quad+\sum_{n=0}^{N-1}\frac{(1-aq^{2n})(b,c,d,e,aq/bg,aq/cg,aq/dg,aq/eg;q)_n}{(aq/b,aq/c,aq/d,aq/e,bf/a,cf/a,df/a,ef/a;q)_n}\frac{(f;q)_{2n}}{(aq/g;q)_{2n}}\left(\frac{fg}{a}\right)^n.\label{eq:finite-sum-q}
\end{align}
\end{Cor}
Under the assumption $|fg/a|<1$, letting $N\to\infty$ in \eqref{eq:finite-sum-q} gives the following corollary.
\begin{Cor}
\label{cor:finite-sum-q-infinity}
Assume, in addition to \eqref{eq:balance}, that $\left|\frac{fg}a\right|<1$.  Then
\begin{align*}
&W_\infty^{(q)}\sum_{n=0}^{\infty}\frac{(1-aq^{2n})(b,c,d,e;q)_n}{(aq/b,aq/c,aq/d,aq/e;q)_n}\left(\frac{fg}{a}\right)^n\\
&\quad=\frac{f}{a}\frac{(1-a/b)(1-a/c)(1-a/d)(1-a/e)(1-a/f)}{1-g/q}\\
&\qquad\cdot\sum_{n=0}^{\infty}\frac{(1-fq^{2n+1}/g)(aq/bg,aq/cg,aq/dg,aq/eg,fq/a,f;q)_n}{(bf/a,cf/a,df/a,ef/a;q)_{n+1}(q^2/g,aq/g;q)_n}q^n\\
&\quad-\frac{f}{a}\sum_{n=0}^{\infty}\frac{(1-fq^{2n+1}/g)(b,c,d,e;q)_{n+1}}{(aq/b,aq/c,aq/d,aq/e;q)_n}\\
&\qquad\cdot\frac{(aq/bg,aq/cg,aq/dg,aq/eg;q)_n}{(bf/a,cf/a,df/a,ef/a;q)_{n+1}}\frac{(f;q)_{2n+1}}{(aq/g;q)_{2n+1}}
\left(\frac{fg}{a}\right)^n\\
&\quad+\sum_{n=0}^{\infty}\frac{(1-aq^{2n})(b,c,d,e,aq/bg,aq/cg,aq/dg,aq/eg;q)_n}{(aq/b,aq/c,aq/d,aq/e,bf/a,cf/a,df/a,ef/a;q)_n}\frac{(f;q)_{2n}}{(aq/g;q)_{2n}}\left(\frac{fg}{a}\right)^n.
\end{align*}
\end{Cor}

The specialization $g=a/f$ lies on the boundary $|fg/a|=1$ and is not covered by Corollary \ref{cor:finite-sum-q-infinity}. We therefore specialize the finite identity first and extract its constant term. To formulate the resulting finite part, set
\begin{align}
H^{(q)}&\coloneqq\frac{(b,c,d,e;q)_\infty}{(aq/b,aq/c,aq/d,aq/e;q)_\infty},\notag\\
C^{(q)}&\coloneqq\frac{(b,c,d,e,fq/b,fq/c,fq/d,fq/e;q)_\infty}{(aq/b,aq/c,aq/d,aq/e,bf/a,cf/a,df/a,ef/a;q)_\infty}.\notag
\end{align}

\begin{Thm}\label{thm:finite-part-g-a-over-f}
Suppose that
\begin{equation}\label{eq:boundary-balance}
 a^2q=bcde.
\end{equation}
Then
\begin{align}
&\frac{(f/a,fq/b,fq/c,fq/d,fq/e;q)_\infty}{(fq/a,bf/a,cf/a,df/a,ef/a;q)_\infty}\notag\\
&\qquad\cdot\Biggl(\sum_{n=0}^{\infty}\left(\frac{(1-aq^{2n})(b,c,d,e;q)_n}{(aq/b,aq/c,aq/d,aq/e;q)_n}-H^{(q)}\right)+H^{(q)}\sum_{j=1}^{\infty}\frac{fq^j/a}{1-fq^j/a}\Biggr)\notag\\
&=-\frac{f^2q}{a^2}(1-a/b)(1-a/c)(1-a/d)(1-a/e)(1-a/f)\notag\\
&\qquad\cdot\sum_{n=0}^{\infty}\frac{(1-f^2q^{2n+1}/a)(fq/b,fq/c,fq/d,fq/e;q)_n}{(bf/a,cf/a,df/a,ef/a;q)_{n+1}(1-fq^n)(1-fq^{n+1}/a)}q^n\notag\\
&\quad-\frac{f}{a}\sum_{n=0}^{\infty}\Biggl(\frac{(1-f^2q^{2n+1}/a)(b,c,d,e;q)_{n+1}}{(aq/b,aq/c,aq/d,aq/e;q)_n}\frac{(fq/b,fq/c,fq/d,fq/e;q)_n}{(bf/a,cf/a,df/a,ef/a;q)_{n+1}}\frac{1}{1-fq^{2n+1}}-C^{(q)}\Biggr)\notag\\
&\quad+\sum_{n=0}^{\infty}\Biggl(\frac{(1-aq^{2n})(b,c,d,e,fq/b,fq/c,fq/d,fq/e;q)_n}{(aq/b,aq/c,aq/d,aq/e,bf/a,cf/a,df/a,ef/a;q)_n}\frac{1}{1-fq^{2n}}-C^{(q)}\Biggr).\label{eq:finite-part-g-a-over-f}
\end{align}
All series in \eqref{eq:finite-part-g-a-over-f} are absolutely convergent.
\end{Thm}

\begin{proof}
Set $g=a/f$ in \eqref{eq:finite-sum-q} and divide by $1-f$ before taking the limit.  The balancing condition \eqref{eq:balance} becomes \eqref{eq:boundary-balance}, and the finite identity takes the form
\begin{align}
&\frac{W_N^{(q)}}{1-f}\sum_{k=0}^{N-1}\frac{(1-aq^{2k})(b,c,d,e;q)_k}{(aq/b,aq/c,aq/d,aq/e;q)_k}\frac{(1-fq^N)(1-fq^N/a)}{(1-fq^{N+k})(1-fq^{N-k}/a)}\notag\\
&\quad=-\frac{f^2q}{a^2}(1-a/b)(1-a/c)(1-a/d)(1-a/e)(1-a/f)\notag\\
&\qquad\cdot\sum_{n=0}^{N-1}\frac{(1-f^2q^{2n+1}/a)(fq/b,fq/c,fq/d,fq/e;q)_n}{(bf/a,cf/a,df/a,ef/a;q)_{n+1}(1-fq^n)(1-fq^{n+1}/a)}q^n
\notag\\
&\qquad-\frac{f}{a}\sum_{n=0}^{N-1}\frac{(1-f^2q^{2n+1}/a)(b,c,d,e;q)_{n+1}}{(aq/b,aq/c,aq/d,aq/e;q)_n}\frac{(fq/b,fq/c,fq/d,fq/e;q)_n}{(bf/a,cf/a,df/a,ef/a;q)_{n+1}}\frac{1}{1-fq^{2n+1}}\notag\\
&\qquad+\sum_{n=0}^{N-1}\frac{(1-aq^{2n})(b,c,d,e,fq/b,fq/c,fq/d,fq/e;q)_n}{(aq/b,aq/c,aq/d,aq/e,bf/a,cf/a,df/a,ef/a;q)_n}\frac{1}{1-fq^{2n}},\label{eq:finite-part-finite-identity}
\end{align}
where we used
\begin{align}
q^k\frac{(fq^N,(a/f)q^{-N};q)_k}{(aq^{1-N}/f,fq^{N+1};q)_k}&=\frac{(1-fq^N)(1-fq^N/a)}{(1-fq^{N+k})(1-fq^{N-k}/a)}.\label{eq:boundary-factor-exact}
\end{align}
Furthermore,
\begin{equation*}
 \lim_{N\to\infty}\frac{W_N^{(q)}}{1-f}H^{(q)}=\frac{(f/a,fq/b,fq/c,fq/d,fq/e;q)_\infty}{(fq/a,bf/a,cf/a,df/a,ef/a;q)_\infty}H^{(q)}=(1-f/a)C^{(q)}.
\end{equation*}

Subtract $N(1-f/a)C^{(q)}$ from both sides of \eqref{eq:finite-part-finite-identity}.  Since $W_N^{(q)}/(1-f)$ approaches its normalized infinite-product limit at rate $O(|q|^N)$, we have $N(W_N^{(q)}-W_\infty^{(q)})/(1-f)\to0$.  Moreover,
\[
 \frac{(1-aq^{2k})(b,c,d,e;q)_k}{(aq/b,aq/c,aq/d,aq/e;q)_k}-H^{(q)}=O(|q|^k),
\]
and \eqref{eq:boundary-factor-exact} shows that its right-hand side is bounded uniformly for $0\leq k<N$.  Hence dominated convergence applies to the part containing the displayed difference.  For the constant part, set $j=N-k$.  Uniformly for $1\leq j\leq N$,
\[
 \frac{(1-fq^N)(1-fq^N/a)}{(1-fq^{2N-j})(1-fq^j/a)}-1=O(|q|^j)+O(|q|^N),
\]
and, for each fixed $j\geq1$,
\[
 \frac{(1-fq^N)(1-fq^N/a)}{(1-fq^{2N-j})(1-fq^j/a)}\longrightarrow\frac{1}{1-fq^j/a}.
\]
These estimates give
\begin{align}
&\lim_{N\to\infty}\Biggl(\frac{W_N^{(q)}}{1-f}\sum_{k=0}^{N-1}\frac{(1-aq^{2k})(b,c,d,e;q)_k}{(aq/b,aq/c,aq/d,aq/e;q)_k}\frac{(1-fq^N)(1-fq^N/a)}{(1-fq^{N+k})(1-fq^{N-k}/a)}-N(1-f/a)C^{(q)}\Biggr)\notag\\
&\quad=\frac{(f/a,fq/b,fq/c,fq/d,fq/e;q)_\infty}{(fq/a,bf/a,cf/a,df/a,ef/a;q)_\infty}\notag\\
&\qquad\cdot\left(\sum_{n=0}^{\infty}\left(\frac{(1-aq^{2n})(b,c,d,e;q)_n}{(aq/b,aq/c,aq/d,aq/e;q)_n}-H^{(q)}\right)+H^{(q)}\sum_{j=1}^{\infty}\left(\frac{1}{1-fq^j/a}-1\right)\right).\label{eq:left-finite-part}
\end{align}

On the right-hand side, the summand in the first series is $O(|q|^n)$, and each of the two expressions
\begin{align*}
&\frac{(1-f^2q^{2n+1}/a)(b,c,d,e;q)_{n+1}}{(aq/b,aq/c,aq/d,aq/e;q)_n}\frac{(fq/b,fq/c,fq/d,fq/e;q)_n}{(bf/a,cf/a,df/a,ef/a;q)_{n+1}}\frac{1}{1-fq^{2n+1}}-C^{(q)},\\
&\frac{(1-aq^{2n})(b,c,d,e;q)_n}{(aq/b,aq/c,aq/d,aq/e;q)_n}\frac{(fq/b,fq/c,fq/d,fq/e;q)_n}{(bf/a,cf/a,df/a,ef/a;q)_n}\frac{1}{1-fq^{2n}}-C^{(q)}
\end{align*}
is $O(|q|^n)$.  Consequently,
\begin{align}
&\lim_{N\to\infty}\Biggl(-\frac{f}{a}\sum_{n=0}^{N-1}\frac{(1-f^2q^{2n+1}/a)(b,c,d,e;q)_{n+1}}{(aq/b,aq/c,aq/d,aq/e;q)_n}\frac{(fq/b,fq/c,fq/d,fq/e;q)_n}{(bf/a,cf/a,df/a,ef/a;q)_{n+1}}\frac{1}{1-fq^{2n+1}}\notag\\
&\hspace{20mm}+\sum_{n=0}^{N-1}\frac{(1-aq^{2n})(b,c,d,e;q)_n}{(aq/b,aq/c,aq/d,aq/e;q)_n}\notag\\
&\hspace{25mm}\cdot\frac{(fq/b,fq/c,fq/d,fq/e;q)_n}{(bf/a,cf/a,df/a,ef/a;q)_n}\frac{1}{1-fq^{2n}}-N(1-f/a)C^{(q)}\Biggr)\notag\\
&\quad=-\frac{f}{a}\sum_{n=0}^{\infty}\Biggl(\frac{(1-f^2q^{2n+1}/a)(b,c,d,e;q)_{n+1}}{(aq/b,aq/c,aq/d,aq/e;q)_n}\notag\\
&\hspace{30mm}\cdot\frac{(fq/b,fq/c,fq/d,fq/e;q)_n}{(bf/a,cf/a,df/a,ef/a;q)_{n+1}}\frac{1}{1-fq^{2n+1}}-C^{(q)}\Biggr)\notag\\
&\qquad+\sum_{n=0}^{\infty}\Biggl(\frac{(1-aq^{2n})(b,c,d,e,fq/b,fq/c,fq/d,fq/e;q)_n}{(aq/b,aq/c,aq/d,aq/e,bf/a,cf/a,df/a,ef/a;q)_n}\frac{1}{1-fq^{2n}}-C^{(q)}\Biggr).\label{eq:right-finite-part}
\end{align}
Finally,
\[
 \frac{1}{1-fq^j/a}-1=\frac{fq^j/a}{1-fq^j/a}.
\]
Letting $N\to\infty$ in \eqref{eq:finite-part-finite-identity} and using \eqref{eq:left-finite-part} and \eqref{eq:right-finite-part} proves \eqref{eq:finite-part-g-a-over-f}.
\end{proof}

\begin{Cor}\label{cor:q-equals-a}
Suppose that
\begin{equation}
 aq=bcd, \label{eq:q-equals-a-balance}
\end{equation}
and that the parameters are such that all expressions below are defined.  Then
\begin{align}
&\sum_{n=0}^{\infty}\Biggl(\frac{(1-aq^{2n})(b,c,d,a,fq/b,fq/c,fq/d,fq/a;q)_n}{(aq/b,aq/c,aq/d,q,bf/a,cf/a,df/a,f;q)_n}\frac{1}{1-fq^{2n}}\notag\\
&\hspace{9mm}-\frac{f}{a}\frac{(1-f^2q^{2n+1}/a)(b,c,d,a;q)_{n+1}}{(aq/b,aq/c,aq/d,q;q)_n}\frac{(fq/b,fq/c,fq/d,fq/a;q)_n}{(bf/a,cf/a,df/a,f;q)_{n+1}}\frac{1}{1-fq^{2n+1}}\notag\\
&\hspace{9mm}-\frac{(b,c,d,a,f/a,fq/b,fq/c,fq/d;q)_\infty}{(aq/b,aq/c,aq/d,q,bf/a,cf/a,df/a,f;q)_\infty}\Biggr)\notag\\
&\quad=\frac{(b,c,d,a,f/a,fq/b,fq/c,fq/d;q)_\infty}{(aq/b,aq/c,aq/d,q,bf/a,cf/a,df/a,f;q)_\infty}\notag\\
&\qquad\cdot\sum_{k=0}^{\infty}\left(\frac{bq^k}{1-bq^k}+\frac{cq^k}{1-cq^k}+\frac{dq^k}{1-dq^k}-\frac{q^{k+1}}{1-q^{k+1}}+\frac{fq^{k+1}/a}{1-fq^{k+1}/a}\right).\label{eq:q-equals-a}
\end{align}
The series in \eqref{eq:q-equals-a} is absolutely convergent.
\end{Cor}

\begin{proof}
Set $e=a$ in Theorem~\ref{thm:finite-part-g-a-over-f}.  The balancing condition becomes \eqref{eq:q-equals-a-balance}, and the first series on the right-hand side of \eqref{eq:finite-part-g-a-over-f} vanishes because its prefactor contains $1-a/e$.  Applying Corollary~\ref{cor:appendix-q} with its parameters $(d,e,f)$ replaced by $(b,c,d)$ gives
\begin{align*}
&\sum_{n=0}^{\infty}\left(\frac{(1-aq^{2n})(b,c,d,a;q)_n}{(aq/b,aq/c,aq/d,q;q)_n}-H^{(q)}\right)\notag\\
&\quad=H^{(q)}\sum_{k=0}^{\infty}\left(\frac{bq^k}{1-bq^k}+\frac{cq^k}{1-cq^k}+\frac{dq^k}{1-dq^k}-\frac{q^{k+1}}{1-q^{k+1}}\right).
\end{align*}
Moreover, after the specialization $e=a$,
\begin{align*}
&\frac{(f/a,fq/b,fq/c,fq/d,fq/e;q)_\infty}{(fq/a,bf/a,cf/a,df/a,ef/a;q)_\infty}H^{(q)}\notag\\
&\quad=\frac{(b,c,d,a,f/a,fq/b,fq/c,fq/d;q)_\infty}{(aq/b,aq/c,aq/d,q,bf/a,cf/a,df/a,f;q)_\infty}.
\end{align*}
Substituting these two formulas into the left-hand side of \eqref{eq:finite-part-g-a-over-f} and combining the second and third series on its right-hand side gives \eqref{eq:q-equals-a}.
\end{proof}

\section{Classical limits}\label{sec:classical-limit}

We finally record the classical limits of Corollaries \ref{cor:first-contiguous-q}--\ref{cor:finite-sum-q-infinity}.  In every formula before the final boundary specialization, the balancing condition is
\begin{equation}
 3a+1=b+c+d+e+f+g. \label{eq:classical-balance}
\end{equation}
For $x\in\CC$ and $n\in\ZZ$, put
\[
 (x)_n\coloneqq \frac{\Gamma(x+n)}{\Gamma(x)},\qquad(x_1,\ldots,x_r)_n\coloneqq \prod_{j=1}^r(x_j)_n,
\]
whenever the expressions involved are defined.  For integers $m\leq M$, we define the classical partial sum by
\begin{equation*}
 \Om_m^M(a\mid b,c,d,e,f,g)\coloneqq \sum_{n=m}^{M-1}\frac{(2n+a)(b,c,d,e,f,g)_n}{(1+a-b,1+a-c,1+a-d,1+a-e,1+a-f,1+a-g)_n}.
\end{equation*}
If the series converges, we use the analogous notation
\begin{equation*}
 \Om_m^\infty(a\mid b,c,d,e,f,g)\coloneqq \sum_{n=m}^{\infty}\frac{(2n+a)(b,c,d,e,f,g)_n}{(1+a-b,1+a-c,1+a-d,1+a-e,1+a-f,1+a-g)_n}.
\end{equation*}
In the unilateral $H$-series notation used in \cite{chu-nonterminating}, the case $m=0$ is
\begin{align}
 &\Om_0^\infty(a\mid b,c,d,e,f,g)\notag\\
 &=a\,{}_7H_7^+\!\left[
 \begin{matrix}1+a/2,b,c,d,e,f,g\\
 a/2,1+a-b,1+a-c,1+a-d,1+a-e,1+a-f,1+a-g
 \end{matrix};1\right].
 \label{eq:classical-omega-7H7}
\end{align}
Thus the classical objects in this section are finite partial sums or unilateral parts of the well-poised ${}_7H_7(1)$ series.  Such unilateral and bilateral series and their Abel-type transformations were studied by Chu \cite[Section~2]{chu-nonterminating}; see also the classical Abel framework in \cite{chu-classical}.

Indeed, after replacing $(a,b,c,d,e,f,g)$ in Section \ref{sec:degeneration} by $(q^a,q^b,q^c,q^d,q^e,q^f,q^g)$, respectively,
we have
\begin{equation*}
 \lim_{q\to1^-}\frac{1}{1-q}\Om_m^M(q^a\mid q^b,q^c,q^d,q^e,q^f,q^g;q)=\Om_m^M(a\mid b,c,d,e,f,g).
\end{equation*}

For the first contiguous relation, set
\begin{align}
 E_r&\coloneqq \frac{e+f+g+r-a}{(e+f-a)(e+g-a)(f+g-a)}\notag\\
 &\qquad\cdot\frac{(b,c,d,e,f,g)_r}{(1+a-b,1+a-c,1+a-d)_r(1+a-e,1+a-f,1+a-g)_{r-1}}.\notag
\end{align}

\begin{Cor}[Classical limit of Corollary \ref{cor:first-contiguous-q}]
\label{cor:first-contiguous-classical} Under \eqref{eq:classical-balance}, for integers $m\leq M$,
\begin{align}
 \Om_m^M(a\mid b,c,d,e,f,g)&=E_m-E_M\notag\\
 &\quad+\frac{efg(1+a-b-c)(1+a-b-d)(1+a-c-d)}{(e+f-a)(e+g-a)(f+g-a)(1+a-b)(1+a-c)(1+a-d)}\notag\\
 &\qquad\cdot\Om_m^M(a+1\mid b,c,d,e+1,f+1,g+1).\notag
\end{align}
\end{Cor}

\begin{Rem}
Take $m=0$ and let $M\to\infty$ in Corollary~\ref{cor:first-contiguous-classical}, retaining the generally nonzero limit of $E_M$. After evaluating this limit using the gamma-ratio asymptotic, applying \eqref{eq:classical-omega-7H7}, and permuting the two parameter triples, we recover the unilateral ${}_7H_7(1)$ recurrence in \cite[Section~2.1]{chu-nonterminating}. Corollary~\ref{cor:first-contiguous-classical} therefore gives a version of this recurrence on an arbitrary finite interval.
\end{Rem}

For the second contiguous relation, put
\begin{align}
 F_r&\coloneqq \frac{g-f-1}{(a-b-f)(a-c-f)(a-d-f)(a-e-f)}\notag\\
 &\quad\cdot\frac{(b,c,d,e,f)_r(g)_{r-1}}{(1+a-b,1+a-c,1+a-d,1+a-e,1+a-f)_{r-1}(1+a-g)_r}.\notag
\end{align}

\begin{Cor}[Classical limit of Corollary \ref{cor:second-contiguous-q}]
\label{cor:second-contiguous-classical}
Under \eqref{eq:classical-balance}, for integers $m\leq M$,
\begin{align}
 \Om_m^M(a\mid b,c,d,e,f,g)&=F_m-F_M\notag\\
 &\quad+
 \frac{f(a-f)(1+a-b-g)(1+a-c-g)(1+a-d-g)(1+a-e-g)}{(g-1)(1+a-g)(a-b-f)(a-c-f)(a-d-f)(a-e-f)}\notag\\
 &\qquad\cdot\Om_m^M(a\mid b,c,d,e,f+1,g-1).\notag
\end{align}
\end{Cor}

The classical transformation formulas below are closely related to the ${}_7H_7(1)$ transformations studied in \cite{chu-nonterminating} and to the use of partial well-poised sums by Chu and Zhang to accelerate Dougall's ${}_5F_4$ summation \cite{chu-zhang-dougall}.

For the iterated formulas, define
\begin{equation*}
 W_N\coloneqq \frac{(f,f-a,1+a-b-g,1+a-c-g,1+a-d-g,1+a-e-g)_N}{(1-g,1+a-g,b+f-a,c+f-a,d+f-a,e+f-a)_N}.
\end{equation*}

\begin{Cor}[Classical limit of Corollary \ref{cor:fixed-lower-q}]
\label{cor:fixed-lower-classical}
Let $N$ be a nonnegative integer and let $m\leq M$ be integers.  Under
\eqref{eq:classical-balance},
\begin{align}
 &W_N\Om_m^M(a\mid b,c,d,e,f+N,g-N)\notag\\
 &=\Om_m^M(a\mid b,c,d,e,f,g)\notag\\
 &\quad+\frac{(b,c,d,e,f)_m(g)_{m-1}}{(1+a-b,1+a-c,1+a-d,1+a-e,1+a-f)_{m-1}(1+a-g)_m}\notag\\
 &\qquad\cdot\sum_{n=0}^{N-1}\frac{(f+2n+1-g)(1+a-b-g,1+a-c-g,1+a-d-g,1+a-e-g)_n}{(b+f-a,c+f-a,d+f-a,e+f-a)_{n+1}}\notag\\
 &\hspace{42mm}\cdot\frac{(f+1-m-a,f+m)_n}{(2-m-g,1+a+m-g)_n}\notag\\
 &\quad-\frac{(b,c,d,e,f)_M(g)_{M-1}}{(1+a-b,1+a-c,1+a-d,1+a-e,1+a-f)_{M-1}(1+a-g)_M}\notag\\
 &\qquad\cdot\sum_{n=0}^{N-1}\frac{(f+2n+1-g)(1+a-b-g,1+a-c-g,1+a-d-g,1+a-e-g)_n}{(b+f-a,c+f-a,d+f-a,e+f-a)_{n+1}}\notag\\
 &\hspace{42mm}\cdot\frac{(f+1-M-a,f+M)_n}{(2-M-g,1+a+M-g)_n}.\notag
\end{align}
\end{Cor}

For the first nonterminating limit, put
\begin{equation*}
 W_\infty\coloneqq \frac{\Gamma(1-g)\Gamma(1+a-g)\Gamma(b+f-a)\Gamma(c+f-a)\Gamma(d+f-a)\Gamma(e+f-a)}{\Gamma(f)\Gamma(f-a)\Gamma(1+a-b-g)\Gamma(1+a-c-g)\Gamma(1+a-d-g)\Gamma(1+a-e-g)}.
\end{equation*}

\begin{Cor}[Classical limit of Corollary \ref{cor:fixed-lower-q-infinity}]\label{cor:fixed-lower-classical-infinity}
Suppose that \eqref{eq:classical-balance} holds, $m\in\ZZ$, and $\operatorname{Re}(f+g-a)>0$.  Then
\begin{align}
 &W_\infty\sum_{k=m}^{\infty}\frac{(2k+a)(b,c,d,e)_k}{(1+a-b,1+a-c,1+a-d,1+a-e)_k}\notag\\
 &=\Om_m^\infty(a\mid b,c,d,e,f,g)\notag\\
 &\quad+\frac{(b,c,d,e,f)_m(g)_{m-1}}{(1+a-b,1+a-c,1+a-d,1+a-e,1+a-f)_{m-1}(1+a-g)_m}\notag\\
 &\qquad\cdot\sum_{n=0}^{\infty}\frac{(f+2n+1-g)(1+a-b-g,1+a-c-g,1+a-d-g,1+a-e-g)_n}{(b+f-a,c+f-a,d+f-a,e+f-a)_{n+1}}\notag\\
 &\hspace{35mm}\cdot\frac{(f+1-m-a,f+m)_n}{(2-m-g,1+a+m-g)_n}\notag\\
 &\quad-\frac{\Gamma(1+a-b)\Gamma(1+a-c)\Gamma(1+a-d)\Gamma(1+a-e)\Gamma(1+a-f)\Gamma(1+a-g)}{\Gamma(b)\Gamma(c)\Gamma(d)\Gamma(e)\Gamma(f)\Gamma(g)}\notag\\
 &\qquad\cdot\sum_{n=0}^{\infty}\frac{(f+2n+1-g)(1+a-b-g,1+a-c-g,1+a-d-g,1+a-e-g)_n}{(b+f-a,c+f-a,d+f-a,e+f-a)_{n+1}}.\notag
\end{align}
\end{Cor}

\begin{Cor}[Classical limit of Corollary \ref{cor:moving-lower-q}]\label{cor:moving-lower-classical}
Let $N$ and $M$ be integers with $0\leq N\leq M$.  Under \eqref{eq:classical-balance},
\begin{align}
 &W_N\Om_N^M(a\mid b,c,d,e,f+N,g-N)\notag\\
 &=\Om_0^M(a\mid b,c,d,e,f,g)\notag\\
 &\quad-\frac{(b,c,d,e,f)_M(g)_{M-1}}{(1+a-b,1+a-c,1+a-d,1+a-e,1+a-f)_{M-1}(1+a-g)_M}\notag\\
 &\qquad\cdot\sum_{n=0}^{N-1}\frac{(f+2n+1-g)(1+a-b-g,1+a-c-g,1+a-d-g,1+a-e-g)_n}{(b+f-a,c+f-a,d+f-a,e+f-a)_{n+1}}\notag\\
 &\hspace{35mm}\cdot\frac{(f+1-M-a,f+M)_n}{(2-M-g,1+a+M-g)_n}\notag\\
 &\quad+\sum_{n=0}^{N-1}\frac{(f+2n+1-g)(b,c,d,e)_{n+1}}{(1+a-b,1+a-c,1+a-d,1+a-e)_n}\notag\\
 &\hspace{28mm}\cdot\frac{(1+a-b-g,1+a-c-g,1+a-d-g,1+a-e-g)_n}{(b+f-a,c+f-a,d+f-a,e+f-a)_{n+1}}\frac{(f)_{2n+1}}{(1+a-g)_{2n+1}}\notag\\
 &\quad-\sum_{n=0}^{N-1}\frac{(2n+a)(b,c,d,e)_n}{(1+a-b,1+a-c,1+a-d,1+a-e)_n}\notag\\
 &\hspace{28mm}\cdot\frac{(1+a-b-g,1+a-c-g,1+a-d-g,1+a-e-g)_n}{(b+f-a,c+f-a,d+f-a,e+f-a)_n}\frac{(f)_{2n}}{(1+a-g)_{2n}}.\notag
\end{align}
\end{Cor}

\begin{Cor}[Classical limit of Corollary \ref{cor:finite-sum-q}]\label{cor:finite-sum-classical}
Let $N$ be a nonnegative integer.  Under \eqref{eq:classical-balance},
\begin{align}
 &W_N\Om_0^N(a\mid b,c,d,e,f+N,g-N)\notag\\
 &=\frac{(a-b)(a-c)(a-d)(a-e)(a-f)}{g-1}\notag\\
 &\qquad\cdot\sum_{n=0}^{N-1}\frac{(f+2n+1-g)(1+a-b-g,1+a-c-g,1+a-d-g,1+a-e-g)_n}{(b+f-a,c+f-a,d+f-a,e+f-a)_{n+1}}\notag\\
 &\hspace{35mm}\cdot\frac{(f+1-a,f)_n}{(2-g,1+a-g)_n}\notag\\
 &\quad-\sum_{n=0}^{N-1}\frac{(f+2n+1-g)(b,c,d,e)_{n+1}}{(1+a-b,1+a-c,1+a-d,1+a-e)_n}\notag\\
 &\hspace{28mm}\cdot\frac{(1+a-b-g,1+a-c-g,1+a-d-g,1+a-e-g)_n}{(b+f-a,c+f-a,d+f-a,e+f-a)_{n+1}}\frac{(f)_{2n+1}}{(1+a-g)_{2n+1}}\notag\\
 &\quad+\sum_{n=0}^{N-1}\frac{(2n+a)(b,c,d,e)_n}{(1+a-b,1+a-c,1+a-d,1+a-e)_n}\notag\\
 &\hspace{28mm}\cdot\frac{(1+a-b-g,1+a-c-g,1+a-d-g,1+a-e-g)_n}{(b+f-a,c+f-a,d+f-a,e+f-a)_n}\frac{(f)_{2n}}{(1+a-g)_{2n}}.\label{eq:finite-sum-classical}
\end{align}
\end{Cor}

\begin{Cor}[Classical limit of Corollary \ref{cor:finite-sum-q-infinity}]\label{cor:finite-sum-classical-infinity}
Suppose that \eqref{eq:classical-balance} holds and $\operatorname{Re}(f+g-a)>0$.  Then
\begin{align}
 &W_\infty\sum_{n=0}^{\infty}\frac{(2n+a)(b,c,d,e)_n}{(1+a-b,1+a-c,1+a-d,1+a-e)_n}\notag\\
 &=\frac{(a-b)(a-c)(a-d)(a-e)(a-f)}{g-1}\notag\\
 &\qquad\cdot\sum_{n=0}^{\infty}\frac{(f+2n+1-g)(1+a-b-g,1+a-c-g,1+a-d-g,1+a-e-g)_n}{(b+f-a,c+f-a,d+f-a,e+f-a)_{n+1}}\notag\\
 &\hspace{35mm}\cdot\frac{(f+1-a,f)_n}{(2-g,1+a-g)_n}\notag\\
&\quad+\sum_{n=0}^{\infty}\Biggl(\frac{(2n+a)(b,c,d,e)_n}{(1+a-b,1+a-c,1+a-d,1+a-e)_n}\notag\\
&\hspace{28mm}\cdot\frac{(1+a-b-g,1+a-c-g,1+a-d-g,1+a-e-g)_n}{(b+f-a,c+f-a,d+f-a,e+f-a)_n}\frac{(f)_{2n}}{(1+a-g)_{2n}}\notag\\
&\hspace{15mm}-\frac{(f+2n+1-g)(b,c,d,e)_{n+1}}{(1+a-b,1+a-c,1+a-d,1+a-e)_n}\notag\\
&\hspace{28mm}\cdot\frac{(1+a-b-g,1+a-c-g,1+a-d-g,1+a-e-g)_n}{(b+f-a,c+f-a,d+f-a,e+f-a)_{n+1}}\frac{(f)_{2n+1}}{(1+a-g)_{2n+1}}\Biggr).\notag
\end{align}
All series displayed above are absolutely convergent.
\end{Cor}

We conclude with the boundary case of Corollary \ref{cor:finite-sum-classical}.  Asymptotic limits of Abel recurrences for Dougall-type bilateral hypergeometric sums were used by Chu \cite{chu-asymptotic}; related transformations of partial well-poised sums were developed by Chu and Zhang \cite{chu-zhang-dougall}.  The argument below uses a similar method to extract
the logarithmic finite part in the boundary case of the present unilateral ${}_7H_7(1)$ series.  Let $\gamma$ denote Euler's constant and let $\psi(x)$ be the digamma function. Define
\begin{align*}
 H&\coloneqq 2\frac{\Gamma(1+a-b)\Gamma(1+a-c)\Gamma(1+a-d)\Gamma(1+a-e)}{\Gamma(b)\Gamma(c)\Gamma(d)\Gamma(e)},\\
 D&\coloneqq 3(f-a)\frac{\Gamma(1+a-b)\Gamma(1+a-c)\Gamma(1+a-d)\Gamma(1+a-e)}{\Gamma(b)\Gamma(c)\Gamma(d)\Gamma(e)}\\
 &\hspace{12mm}\cdot\frac{\Gamma(b+f-a)\Gamma(c+f-a)\Gamma(d+f-a)\Gamma(e+f-a)}{\Gamma(1+f-b)\Gamma(1+f-c)\Gamma(1+f-d)\Gamma(1+f-e)}.
\end{align*}

\begin{Thm}\label{thm:classical-finite-part-g-a-minus-f}
Suppose that
\begin{equation}
 b+c+d+e=2a+1,\label{eq:classical-boundary-balance}
\end{equation}
and that the parameters are such that all expressions below are defined. Then
\begin{align}
 &\frac{\Gamma(b+f-a)\Gamma(c+f-a)\Gamma(d+f-a)\Gamma(e+f-a)}{\Gamma(1+f-b)\Gamma(1+f-c)\Gamma(1+f-d)\Gamma(1+f-e)}\notag\\
 &\qquad\cdot\sum_{n=0}^{\infty}\Biggl(\frac{(2n+a)(b,c,d,e)_n}{(1+a-b,1+a-c,1+a-d,1+a-e)_n}-\frac{H}{n+1}\Biggr)\notag\\
 &\quad-\frac{\Gamma(b+f-a)\Gamma(c+f-a)\Gamma(d+f-a)\Gamma(e+f-a)}{\Gamma(1+f-b)\Gamma(1+f-c)\Gamma(1+f-d)\Gamma(1+f-e)}\notag\\
 &\hspace{12mm}\cdot\frac{\Gamma(1+a-b)\Gamma(1+a-c)\Gamma(1+a-d)\Gamma(1+a-e)}{\Gamma(b)\Gamma(c)\Gamma(d)\Gamma(e)}\bigl(\gamma+\log 2+\psi(1+f-a)\bigr)\notag\\
 &=(a-b)(a-c)(a-d)(a-e)\notag\\
 &\qquad\cdot\sum_{n=0}^{\infty}\frac{(2n+1+2f-a)(1+f-b,1+f-c,1+f-d,1+f-e)_n}{(b+f-a,c+f-a,d+f-a,e+f-a)_{n+1}}\notag\\
 &\hspace{35mm}\cdot\frac{1}{(n+f)(n+1+f-a)}\notag\\
 &\qquad+\frac{1}{f-a}\sum_{n=0}^{\infty}\Biggl(\frac{(2n+a)(b,c,d,e)_n}{(1+a-b,1+a-c,1+a-d,1+a-e)_n}\notag\\
 &\hspace{28mm}\cdot\frac{(1+f-b,1+f-c,1+f-d,1+f-e)_n}{(b+f-a,c+f-a,d+f-a,e+f-a)_n}\frac{1}{2n+f}\notag\\
 &\hspace{15mm}-\frac{(2n+1+2f-a)(b,c,d,e)_{n+1}}{(1+a-b,1+a-c,1+a-d,1+a-e)_n}\notag\\
 &\hspace{28mm}\cdot\frac{(1+f-b,1+f-c,1+f-d,1+f-e)_n}{(b+f-a,c+f-a,d+f-a,e+f-a)_{n+1}}\frac{1}{2n+1+f}-\frac{D}{n+1}\Biggr).\label{eq:classical-finite-part-g-a-minus-f}
\end{align}
All three series in \eqref{eq:classical-finite-part-g-a-minus-f} are absolutely convergent.
\end{Thm}

\begin{proof}
Set $g=a-f$ in Corollary \ref{cor:finite-sum-classical}, divide the resulting identity by $f(f-a)$, and then let $N\to\infty$.  The balancing condition \eqref{eq:classical-balance} becomes \eqref{eq:classical-boundary-balance}.  Throughout the proof we use
$(f)_m/(1+f)_m=f/(m+f)$.  After division, the left-hand side of \eqref{eq:finite-sum-classical} becomes
\begin{align*}
 &\frac{W_N}{f(f-a)}\sum_{k=0}^{N-1}\frac{(2k+a)(b,c,d,e)_k}{(1+a-b,1+a-c,1+a-d,1+a-e)_k}\frac{(N+f)(N+f-a)}{(N+f+k)(N+f-a-k)}.
\end{align*}
Here we used the elementary identity $(x)_k/(x+1)_k=x/(x+k)$ in the form
\begin{equation}
 \frac{(f+N,a-f-N)_k}{(1+a-f-N,1+f+N)_k}=\frac{(N+f)(N+f-a)}{(N+f+k)(N+f-a-k)}.\label{eq:classical-boundary-factor}
\end{equation}
Moreover, the gamma-ratio asymptotic yields
\begin{align*}
 \frac{W_N}{f(f-a)}&=\frac{\Gamma(b+f-a)\Gamma(c+f-a)\Gamma(d+f-a)\Gamma(e+f-a)}{\Gamma(1+f-b)\Gamma(1+f-c)\Gamma(1+f-d)\Gamma(1+f-e)}+O(N^{-1}).
\end{align*}

By \eqref{eq:classical-boundary-balance},
\begin{equation*}
 \frac{(2k+a)(b,c,d,e)_k}{(1+a-b,1+a-c,1+a-d,1+a-e)_k}=\frac{H}{k}+O(k^{-2})\qquad(k\to\infty).
\end{equation*}
From \eqref{eq:classical-boundary-factor},
\[
 \frac{(N+f)(N+f-a)}{(N+f+k)(N+f-a-k)}=\frac{(N+f)(N+f-a)}{2N+2f-a}\left(\frac{1}{N+f+k}+\frac{1}{N+f-a-k}\right).
\]
Partial fractions therefore give
\begin{align}
 &\sum_{k=0}^{N-1}\frac{1}{k+1}\frac{(N+f)(N+f-a)}{(N+f+k)(N+f-a-k)}\notag\\
 &\quad=\frac{(N+f)(N+f-a)}{2N+2f-a}\left(\frac{1}{N+f-1}+\frac{1}{N+f-a+1}\right)\sum_{j=1}^N\frac{1}{j}\notag\\
 &\quad-\frac{(N+f)(N+f-a)}{(2N+2f-a)(N+f-1)}\bigl(\psi(2N+f)-\psi(N+f)\bigr)\notag\\
 &\quad+\frac{(N+f)(N+f-a)}{(2N+2f-a)(N+f-a+1)}\bigl(\psi(N+f-a+1)-\psi(1+f-a)\bigr).\notag
\end{align}
The standard asymptotics of the harmonic numbers and the digamma function imply
\begin{align}
 \sum_{k=0}^{N-1}\frac{1}{k+1}\frac{(N+f)(N+f-a)}{(N+f+k)(N+f-a-k)}
 &=\frac32\log N+\gamma-\frac12\bigl(\log2+\psi(1+f-a)\bigr)+o(1).\label{eq:classical-boundary-harmonic-asymptotic}
\end{align}
On writing
\begin{align*}
&\frac{(2k+a)(b,c,d,e)_k}{(1+a-b,1+a-c,1+a-d,1+a-e)_k}\\
&\quad=\frac{H}{k+1}+\left(\frac{(2k+a)(b,c,d,e)_k}{(1+a-b,1+a-c,1+a-d,1+a-e)_k}-\frac{H}{k+1}\right),
\end{align*}
the expression in parentheses is $O((k+1)^{-2})$.  If the factor in \eqref{eq:classical-boundary-factor} is denoted by $B_{N,k}$, then, with a constant independent of $N$,
\[
 |B_{N,k}|\leq C\quad\left(0\leq k\leq \frac N2\right),\qquad |B_{N,N-j}|\leq \frac{CN}{j+1}\quad\left(1\leq j<\frac N2\right).
\]
The contribution from the second range is therefore $O(N^{-1}\log N)$, while the first range is handled by the summable majorant
$C(k+1)^{-2}$.  It follows that
\begin{align*}
 &\sum_{k=0}^{N-1}\left(\frac{(2k+a)(b,c,d,e)_k}{(1+a-b,1+a-c,1+a-d,1+a-e)_k}-\frac{H}{k+1}\right)\frac{(N+f)(N+f-a)}{(N+f+k)(N+f-a-k)}\notag\\
 &\qquad\longrightarrow\sum_{k=0}^{\infty}\left(\frac{(2k+a)(b,c,d,e)_k}{(1+a-b,1+a-c,1+a-d,1+a-e)_k}-\frac{H}{k+1}\right).
\end{align*}
Equations \eqref{eq:classical-boundary-factor} and \eqref{eq:classical-boundary-harmonic-asymptotic} then give
\begin{align}
 &\frac{W_N}{f(f-a)}\sum_{k=0}^{N-1}\frac{(2k+a)(b,c,d,e)_k}{(1+a-b,1+a-c,1+a-d,1+a-e)_k}\frac{(N+f)(N+f-a)}{(N+f+k)(N+f-a-k)}\notag\\
 &\quad=\frac{D}{f-a}\log N+\frac{\Gamma(b+f-a)\Gamma(c+f-a)\Gamma(d+f-a)\Gamma(e+f-a)}{\Gamma(1+f-b)\Gamma(1+f-c)\Gamma(1+f-d)\Gamma(1+f-e)}\Biggl(\notag\\
 &\qquad\sum_{k=0}^{\infty}\Biggl(\frac{(2k+a)(b,c,d,e)_k}{(1+a-b,1+a-c,1+a-d,1+a-e)_k}-\frac{H}{k+1}\Biggr)\notag\\
 &\hspace{36mm}+H\gamma-\frac{H}{2}\bigl(\log2+\psi(1+f-a)\bigr)\Biggr)+o(1).\label{eq:classical-boundary-left-asymptotic}
\end{align}

It remains to evaluate the right-hand side.  The gamma-ratio asymptotic gives
\begin{align*}
&\frac{(2n+1+2f-a)(1+f-b,1+f-c,1+f-d,1+f-e)_n}{(b+f-a,c+f-a,d+f-a,e+f-a)_{n+1}}\notag\\
&\hspace{18mm}\cdot\frac{1}{(n+f)(n+1+f-a)}=O(n^{-3}),\\
&\frac{(2n+a)(b,c,d,e)_n}{(1+a-b,1+a-c,1+a-d,1+a-e)_n}\notag\\
&\hspace{18mm}\cdot\frac{(1+f-b,1+f-c,1+f-d,1+f-e)_n}{(b+f-a,c+f-a,d+f-a,e+f-a)_n}\frac{1}{2n+f}\longrightarrow \frac{D}{3(f-a)}.
\end{align*}
Furthermore, the ratio of the summand in the second correction sum to that in the third correction sum in Corollary \ref{cor:finite-sum-classical} is
\begin{align}
 \frac{2n+1+2f-a}{2n+a}\frac{2n+f}{2n+f+1}\prod_{x\in(b,c,d,e)}\frac{n+x}{n+x+f-a}&=1+\frac{3(a-f)}{n}+O(n^{-2}).\notag
\end{align}
By the definition of $D$, it follows that
\begin{align*}
&\frac{(2n+a)(b,c,d,e)_n}{(1+a-b,1+a-c,1+a-d,1+a-e)_n}\notag\\
&\hspace{18mm}\cdot\frac{(1+f-b,1+f-c,1+f-d,1+f-e)_n}{(b+f-a,c+f-a,d+f-a,e+f-a)_n}\frac{1}{2n+f}\\
&\quad-\frac{(2n+1+2f-a)(b,c,d,e)_{n+1}}{(1+a-b,1+a-c,1+a-d,1+a-e)_n}\\
&\hspace{18mm}\cdot\frac{(1+f-b,1+f-c,1+f-d,1+f-e)_n}{(b+f-a,c+f-a,d+f-a,e+f-a)_{n+1}}\frac{1}{2n+1+f}\\
&\qquad=\frac{D}{n}+O(n^{-2}).
\end{align*}
Consequently, the right-hand side of \eqref{eq:finite-sum-classical}, specialized to $g=a-f$ and divided by $f(f-a)$, is
\begin{align}
 &\frac{D}{f-a}\log N+(a-b)(a-c)(a-d)(a-e)\notag\\
 &\qquad\cdot\sum_{n=0}^{\infty}\frac{(2n+1+2f-a)(1+f-b,1+f-c,1+f-d,1+f-e)_n}{(b+f-a,c+f-a,d+f-a,e+f-a)_{n+1}}\notag\\
 &\hspace{35mm}\cdot\frac{1}{(n+f)(n+1+f-a)}\notag\\
 &\quad+\frac{1}{f-a}\sum_{n=0}^{\infty}\Biggl(\frac{(2n+a)(b,c,d,e)_n}{(1+a-b,1+a-c,1+a-d,1+a-e)_n}\notag\\
 &\hspace{25mm}\cdot\frac{(1+f-b,1+f-c,1+f-d,1+f-e)_n}{(b+f-a,c+f-a,d+f-a,e+f-a)_n}\frac{1}{2n+f}\notag\\
 &\hspace{12mm}-\frac{(2n+1+2f-a)(b,c,d,e)_{n+1}}{(1+a-b,1+a-c,1+a-d,1+a-e)_n}\notag\\
 &\hspace{25mm}\cdot\frac{(1+f-b,1+f-c,1+f-d,1+f-e)_n}{(b+f-a,c+f-a,d+f-a,e+f-a)_{n+1}}\frac{1}{2n+1+f}-\frac{D}{n+1}\Biggr)\notag\\
 &\quad+\frac{D}{f-a}\gamma+o(1).\label{eq:classical-boundary-right-asymptotic}
\end{align}
Comparing the constant terms in \eqref{eq:classical-boundary-left-asymptotic} and \eqref{eq:classical-boundary-right-asymptotic} using the definition of $D$ proves \eqref{eq:classical-finite-part-g-a-minus-f}.  The estimates above also prove the asserted absolute convergence.
\end{proof}

\begin{Cor}\label{cor:classical-equals-a}
Suppose that
\begin{equation}
 b+c+d=a+1, \label{eq:classical-equals-a-balance}
\end{equation}
and that the parameters are such that all expressions below are defined.  Then
\begin{align}
&\sum_{n=0}^{\infty}\Biggl(\frac{(2n+a)(a,b,c,d,1+f-a,1+f-b,1+f-c,1+f-d)_n}{n!(1+a-b,1+a-c,1+a-d,f,b+f-a,c+f-a,d+f-a)_n}\frac{1}{2n+f}\notag\\
&\hspace{9mm}-\frac{(2n+1+2f-a)(a,b,c,d)_{n+1}}{n!(1+a-b,1+a-c,1+a-d)_n}\notag\\
&\hspace{21mm}\cdot\frac{(1+f-a,1+f-b,1+f-c,1+f-d)_n}{(f,b+f-a,c+f-a,d+f-a)_{n+1}}\frac{1}{2n+1+f}\notag\\
&\hspace{9mm}-\frac{3(f-a)}{n+1} \frac{\Gamma(1+a-b)\Gamma(1+a-c)\Gamma(1+a-d)}{\Gamma(a)\Gamma(b)\Gamma(c)\Gamma(d)}\notag\\
&\hspace{30mm}\cdot\frac{\Gamma(f)\Gamma(b+f-a)\Gamma(c+f-a)\Gamma(d+f-a)}{\Gamma(1+f-a)\Gamma(1+f-b)\Gamma(1+f-c)\Gamma(1+f-d)}\Biggr)
\notag\\
&\quad=-\frac{\Gamma(1+a-b)\Gamma(1+a-c)\Gamma(1+a-d)\Gamma(f)\Gamma(b+f-a)\Gamma(c+f-a)\Gamma(d+f-a)}{\Gamma(a)\Gamma(b)\Gamma(c)\Gamma(d)\Gamma(f-a)\Gamma(1+f-b)\Gamma(1+f-c)\Gamma(1+f-d)}\notag\\
&\hspace{18mm}\cdot\bigl(4\gamma+\log2+\psi(b)+\psi(c)+\psi(d)+\psi(1+f-a)\bigr).\label{eq:classical-equals-a}
\end{align}
The series in \eqref{eq:classical-equals-a} is absolutely convergent.
\end{Cor}

\begin{proof}
Set $e=a$ in Theorem~\ref{thm:classical-finite-part-g-a-minus-f}.  The balancing condition becomes \eqref{eq:classical-equals-a-balance}.  Applying Corollary~\ref{cor:appendix-classical} with its parameters $(d,e,f)$ replaced by $(b,c,d)$ gives
\begin{align*}
&\sum_{n=0}^{\infty}\left(\frac{(2n+a)(a,b,c,d)_n}{(1,1+a-b,1+a-c,1+a-d)_n}-\frac{H}{n+1}\right)\notag\\
&\quad= \frac{\Gamma(1+a-b)\Gamma(1+a-c)\Gamma(1+a-d)}{\Gamma(a)\Gamma(b)\Gamma(c)\Gamma(d)}\bigl(3\psi(1)-\psi(b)-\psi(c)-\psi(d)\bigr).
\end{align*}
The first series on the right-hand side of \eqref{eq:classical-finite-part-g-a-minus-f} vanishes because its prefactor contains $a-e$.  Substitute the preceding evaluation into the left-hand side, use $\psi(1)=-\gamma$, and simplify the defining gamma product for $D$ at $e=a$.  Rearranging the resulting identity gives \eqref{eq:classical-equals-a}.
\end{proof}

\appendix

\section{Well-poised Evans--Stanton-type formulas}\label{sec:appendix-evans-stanton}
In this appendix, we establish well-poised analogues of the zero-balanced finite-part identities studied by Evans and Stanton \cite[Theorems~1 and~2]{evans-stanton}.

\begin{Thm}\label{thm:appendix-first-q}
Suppose that
\[
 a^2q^2=bcdef,\qquad w=\frac{a^2q}{bcd},
\]
and $\left|\frac aw\right|<1$.  Then
\begin{align}
&\sum_{n=0}^{\infty}\left(\frac{(1-aq^{2n})(aq^{n+1}/b,aq^{n+1}/c,aq^{n+1}/d,aq^{n+1}/e,aq^{n+1}/f,q^{n+1};q)_{\infty}}{(aq^n,bq^n,cq^n,dq^n,eq^n,fq^n;q)_{\infty}}-1\right)\notag\\
&\quad=-\frac{(aq/b,aq/c,aq/d,a/w;q)_{\infty}}{(b,c,d,wq;q)_{\infty}}\sum_{n=0}^{\infty}\frac{1-wq^{2n}}{1-w}\frac{(w,wb/a,wc/a,wd/a;q)_n}{(aq/b,aq/c,aq/d,q;q)_n}\left(\frac{a}{w}\right)^n\notag\\
&\hspace{42mm}\cdot\sum_{k=0}^{n-1}\frac{1-wq^{2k+1}}{(1-eq^k)(1-fq^k)}\notag\\
&\qquad+\sum_{k=0}^{\infty}\left(\frac{aq^k/w}{1-aq^k/w}+\frac{eq^k}{1-eq^k}+\frac{fq^k}{1-fq^k}-\frac{q^{k+1}}{1-q^{k+1}}\right).\label{eq:appendix-first-q}
\end{align}
\end{Thm}

\begin{proof}
We begin with Bailey's nonterminating ${}_8\phi_7$ transformation \cite[Eq.~(III.23)]{gr}:
\begin{align*}
&\sum_{n=0}^{\infty}\frac{(1-aq^{2n})(a,b,c,d,e,f;q)_n}{(1-a)(aq/b,aq/c,aq/d,aq/e,aq/f,q;q)_n}\left(\frac{a^2q^2}{bcdef}\right)^n\\
&=\frac{(aq,aq/ef,wq/e,wq/f;q)_{\infty}}{(aq/e,aq/f,wq,wq/ef;q)_{\infty}}\sum_{n=0}^{\infty}\frac{(1-wq^{2n})(w,wb/a,wc/a,wd/a,e,f;q)_n}{(1-w)(aq/b,aq/c,aq/d,wq/e,wq/f,q;q)_n}\left(\frac{aq}{ef}\right)^n.
\end{align*}
Keep $a,b,c,d,e$ fixed and regard $f$ as a function of $t$, where
\[
 t=\frac{a^2q^2}{bcdef}=\frac{wq}{ef}.
\]
Subtracting the limiting coefficient times $\sum_{n=0}^{\infty}t^n=\frac 1{1-t}$ gives
\begin{align}
&\sum_{n=0}^{\infty}\left(\frac{(1-aq^{2n})(a,b,c,d,e,f;q)_n}{(aq/b,aq/c,aq/d,aq/e,aq/f,q;q)_n}-\frac{(a,b,c,d,e,f;q)_{\infty}}{(aq/b,aq/c,aq/d,aq/e,aq/f,q;q)_{\infty}}\right)t^n\notag\\
&=\frac{(a,aq/ef,wq/e,wq/f;q)_{\infty}}{(aq/e,aq/f,wq,t;q)_{\infty}}\sum_{n=0}^{\infty}\frac{(1-wq^{2n})(w,wb/a,wc/a,wd/a,e,f;q)_n}{(1-w)(aq/b,aq/c,aq/d,wq/e,wq/f,q;q)_n}\left(\frac{aq}{ef}\right)^n\notag\\
&\qquad-\frac{(a,b,c,d,e,f;q)_{\infty}}{(aq/b,aq/c,aq/d,aq/e,aq/f,q;q)_{\infty}(1-t)}.
\label{eq:appendix-regularized-series}
\end{align}
Rogers's ${}_6\phi_5$ summation \cite[Eq.~(II.20)]{gr} gives
\[
\sum_{n=0}^{\infty}\frac{(1-wq^{2n})(w,wb/a,wc/a,wd/a;q)_n}{(1-w)(aq/b,aq/c,aq/d,q;q)_n}\left(\frac aw\right)^n=\frac{(wq,b,c,d;q)_{\infty}}{(aq/b,aq/c,aq/d,a/w;q)_{\infty}}.
\]
Using this evaluation, we can rewrite the right-hand side
of \eqref{eq:appendix-regularized-series} as
\begin{align*}
&\frac{(a,aq/ef,wq/e,wq/f;q)_{\infty}}{(aq/e,aq/f,wq,t;q)_{\infty}}\sum_{n=0}^{\infty}\frac{1-wq^{2n}}{1-w}\frac{(w,wb/a,wc/a,wd/a;q)_n}{(aq/b,aq/c,aq/d,q;q)_n}\left(\frac aw\right)^n\left(\frac{(e,wq/et;q)_nt^n}{(wq/e,et;q)_n}-1\right)\\
&\qquad+\frac{(a,b,c,d;q)_{\infty}}{(aq/b,aq/c,aq/d,aq/e,aq/f;q)_{\infty}(1-t)}\left(\frac{(at/w,wq/e,et;q)_{\infty}}{(tq,a/w;q)_{\infty}}-\frac{(e,wq/et;q)_{\infty}}{(q;q)_{\infty}}\right).
\end{align*}
Since $ef=wq$ at $t=1$, logarithmic differentiation gives
\begin{align*}
\left.\frac{d}{dt}\frac{(e,wq/et;q)_n t^n}{(wq/e,et;q)_n}\right|_{t=1}
&=\sum_{k=0}^{n-1}\left(1+\frac{eq^k}{1-eq^k}+\frac{fq^k}{1-fq^k}\right)\notag\\
&=\sum_{k=0}^{n-1}\frac{1-wq^{2k+1}}{(1-eq^k)(1-fq^k)}.
\end{align*}
For the infinite products, we use
\[
 \left.\frac{d}{dt}\log(xt;q)_{\infty}\right|_{t=1}=-\sum_{k=0}^{\infty}\frac{xq^k}{1-xq^k}.
\]
Thus, letting $t\to1$, we obtain
\begin{align*}
&\sum_{n=0}^{\infty}\left(\frac{(1-aq^{2n})(a,b,c,d,e,f;q)_n}{(aq/b,aq/c,aq/d,aq/e,aq/f,q;q)_n}-\frac{(a,b,c,d,e,f;q)_{\infty}}{(aq/b,aq/c,aq/d,aq/e,aq/f,q;q)_{\infty}}\right)\\
&=-\frac{(a,aq/ef,wq/e,wq/f;q)_{\infty}}{(aq/e,aq/f,wq,q;q)_{\infty}}\\
&\qquad\cdot\sum_{n=0}^{\infty}\frac{1-wq^{2n}}{1-w}\frac{(w,wb/a,wc/a,wd/a;q)_n}{(aq/b,aq/c,aq/d,q;q)_n}\left(\frac aw\right)^n\sum_{k=0}^{n-1}\frac{1-wq^{2k+1}}{(1-eq^k)(1-fq^k)}\\
&\qquad-\frac{(a,b,c,d,e,f;q)_{\infty}}{(aq/b,aq/c,aq/d,aq/e,aq/f,q;q)_{\infty}}\\
&\qquad\qquad\cdot\sum_{k=0}^{\infty}\left(\frac{q^{k+1}}{1-q^{k+1}}-\frac{aq^k/w}{1-aq^k/w}-\frac{eq^k}{1-eq^k}-\frac{wq^{k+1}/e}{1-wq^{k+1}/e}\right).
\end{align*}
Dividing both sides by
\[
 \frac{(a,b,c,d,e,f;q)_{\infty}}
 {(aq/b,aq/c,aq/d,aq/e,aq/f,q;q)_{\infty}},
\]
we obtain \eqref{eq:appendix-first-q}.
\end{proof}

Setting $bc=aq$ in Theorem~\ref{thm:appendix-first-q} gives the following corollary.
\begin{Cor}\label{cor:appendix-q}
If $aq=def$, then
\begin{align*}
&\sum_{n=0}^{\infty}\left(\frac{(1-aq^{2n})(aq^{n+1}/d,aq^{n+1}/e,aq^{n+1}/f,q^{n+1};q)_{\infty}}{(aq^n,dq^n,eq^n,fq^n;q)_{\infty}}-1\right)\notag\\
&\quad=\sum_{k=0}^{\infty}\left(\frac{dq^k}{1-dq^k}+\frac{eq^k}{1-eq^k}+\frac{fq^k}{1-fq^k}-\frac{q^{k+1}}{1-q^{k+1}}\right).
\end{align*}
\end{Cor}
To derive the following corollary, we first rewrite the left-hand side of Theorem~\ref{thm:appendix-first-q} as
\begin{align*}
&\sum_{n=0}^{\infty}\left(\frac{(1-aq^{2n})(aq^{n+1}/b,aq^{n+1}/c,aq^{n+1}/d,aq^{n+1}/e,aq^{n+1}/f,q^{n+1};q)_{\infty}}{(aq^n,bq^n,cq^n,dq^n,eq^n,fq^n;q)_{\infty}}-1\right)\\
&=\sum_{n=0}^{\infty}\left(\frac{(1-aq^{2n})(aq^{n+1}/b,aq^{n+1}/c,aq^{n+1}/d,aq^{n+1}/e,aq^{n+1}/f,q^{n+1};q)_{\infty}}{(aq^n,bq^n,cq^n,dq^n,eq^n,fq^n;q)_{\infty}}-\frac{2q^{n+1}}{1-q^{n+1}}-1\right)\\
&\qquad+2\sum_{n=0}^{\infty}\frac{q^{n+1}}{1-q^{n+1}}.
\end{align*}
We then subtract
\[
2\sum_{n=0}^{\infty}\frac{q^{n+1}}{1-q^{n+1}}
\]
from both sides. We next replace $a,b,c,d,e,f$ by $q^a,q^b,q^c,q^d,q^e,q^f$, respectively, multiply both sides by $1-q$, and let $q\to 1$.
\begin{Cor}\label{cor:appendix-first-classical}
Let
\[
 b+c+d+e+f=2a+2,\qquad w=1+2a-b-c-d.
\]
If $\operatorname{Re}(a-w)>0$ and $\operatorname{Re}b>0$, then
\begin{align*}
&\sum_{n=0}^{\infty}\Biggl(\frac{(a+2n)\Gamma(a+n)\Gamma(b+n)\Gamma(c+n)}{\Gamma(1+a-b+n)\Gamma(1+a-c+n)\Gamma(n+1)}\notag\\
&\hspace{18mm}\cdot\frac{\Gamma(d+n)\Gamma(e+n)\Gamma(f+n)}{\Gamma(1+a-d+n)\Gamma(1+a-e+n)\Gamma(1+a-f+n)}-\frac{2}{n+1}\Biggr)
\notag\\
&\quad=-\frac{\Gamma(b)\Gamma(c)\Gamma(d)\Gamma(w)}{\Gamma(1+a-b)\Gamma(1+a-c)\Gamma(1+a-d)\Gamma(a-w)}\notag\\
&\qquad\cdot\sum_{n=0}^{\infty}\frac{(w+2n)(w,w+b-a,w+c-a,w+d-a)_n}{(1+a-b,1+a-c,1+a-d)_n n!}\sum_{k=0}^{n-1}\left(\frac{1}{e+k}+\frac{1}{f+k}\right)\notag\\
&\qquad+\sum_{k=0}^{\infty}\left(\frac{1}{k+a-w}+\frac{1}{e+k}+\frac{1}{f+k}-\frac{3}{k+1}\right)\label{eq:appendix-first-classical}
\end{align*}
\end{Cor}

Applying the same procedure to Corollary \ref{cor:appendix-q} gives the following corollary.
\begin{Cor}\label{cor:appendix-classical}
If $d+e+f=a+1$, then
\begin{align*}
&\sum_{n=0}^{\infty}\left(\frac{(a+2n)\Gamma(a+n)\Gamma(d+n)\Gamma(e+n)\Gamma(f+n)}{\Gamma(n+1)\Gamma(1+a-d+n)\Gamma(1+a-e+n)\Gamma(1+a-f+n)}-\frac{2}{n+1}\right)\notag\\
&\quad=\sum_{k=0}^{\infty}\left(\frac{1}{d+k}+\frac{1}{e+k}+\frac{1}{f+k}-\frac{3}{k+1}\right).\notag
\end{align*}
\end{Cor}


\begin{thebibliography}{99}

\bibitem[C1]{chu-classical}
W. Chu, \emph{Abel's method on summation by parts and hypergeometric series}, J. Difference Equ. Appl. \textbf{12} (2006), no.~8, 783--798.

\bibitem[C2]{chu}
W. Chu, \emph{Abel's lemma on summation by parts and basic hypergeometric series}, Adv. in Appl. Math. \textbf{39} (2007), no.~4, 490--514.

\bibitem[C3]{chu-asymptotic}
W. Chu, \emph{Asymptotic method for Dougall's bilateral hypergeometric sums}, Bull. Sci. Math. \textbf{131} (2007), no.~5, 457--468.

\bibitem[C4]{chu-nonterminating}
W. Chu, \emph{Nonterminating well-poised hypergeometric series}, Hacet. J. Math. Stat. \textbf{49} (2020), no.~5, 1611--1624.

\bibitem[CJ1]{chu-jia}
W. Chu and C. Jia, \emph{Abel's method on summation by parts and theta hypergeometric series}, J. Combin. Theory Ser. A \textbf{115} (2008), no.~5, 815--844.

\bibitem[CJ2]{chu-jia-elliptic}
W. Chu and C. Jia, \emph{Abel's method on summation by parts for elliptic hypergeometric series}, Comm. Contemp. Math. \textbf{11} (2009), no.~3, 337--353.

\bibitem[CW]{chu-wang}
W. Chu and C. Wang, \emph{Abel's lemma on summation by parts and partial $q$-series transformations}, Sci. China Ser. A \textbf{52} (2009), no.~4, 720--748.

\bibitem[CZ1]{chu-zhang-nonterminating}
W. Chu and W. Zhang, \emph{Abel's method on summation by parts and nonterminating $q$-series identities}, Collect. Math. \textbf{60} (2009), no.~2, 193--211.

\bibitem[CZ2]{chu-zhang-bailey}
W. Chu and W. Zhang, \emph{Partial sums of Bailey's ${}_6\psi_6$-series to faster convergent series}, J. Difference Equ. Appl. \textbf{18} (2012), no.~2, 239--260.

\bibitem[CZ3]{chu-zhang-dougall}
W. Chu and W. Zhang, \emph{Accelerating Dougall's ${}_5F_4$-sum and infinite series involving $\pi$}, Math. Comp. \textbf{83} (2014), no.~285, 475--512.

\bibitem[ES]{evans-stanton}
R. J. Evans and D. Stanton, \emph{Asymptotic formulas for zero-balanced hypergeometric series}, SIAM J. Math. Anal. \textbf{15} (1984), no. 5, 1010--1020.

\bibitem[GR]{gr}
G. Gasper and M. Rahman, \emph{Basic Hypergeometric Series}, 2nd ed., Encyclopedia of Mathematics and its Applications, vol.~96, Cambridge University Press, Cambridge, 2004.

\bibitem[S]{schlosser}
M. Schlosser, \emph{Elementary derivations of identities for bilateral basic hypergeometric series}, Selecta Math. (N.S.) \textbf{9} (2003), no.~1, 119--159.

\bibitem[W]{warnaar}
S. O. Warnaar, \emph{Summation and transformation formulas for elliptic hypergeometric series}, Constr. Approx. \textbf{18} (2002), no.~4, 479--502.

\bibitem[WX]{wang-xu}
C. Wang and J. Xu, \emph{New contiguous relations for a balanced $q$-series and their applications}, Ramanujan J. \textbf{64} (2024),  no. 1, 275--287.

\bibitem[XM]{xu-ma}
J. Xu and X. Ma, \emph{General $q$-series transformations based on Abel's lemma on summation by parts and their applications}, J. Difference Equ. Appl. \textbf{30} (2024), no.~5, 553--576.


\end{thebibliography}
\end{document}